\documentclass[11pt,reqno]{amsart}
\usepackage{amsfonts}

\usepackage{fancyhdr}
\usepackage[latin1]{inputenc}
\usepackage{amsmath}
\usepackage{thmdefs}
\usepackage{natbib}
\usepackage{graphicx}
\usepackage{hyperref}
\usepackage{dsfont}

\theoremstyle{definition} \theoremstyle{remark}
\numberwithin{equation}{section}

\renewcommand{\cite}{\citet}

\begin{document}
\title{Parametric estimation and tests through divergences and duality
technique}
\author{Michel Broniatowski$^{*}$ and Amor Keziou$^{**}$}
\address{$^{*}$LSTA-Université Paris 6, e-mail: mbr@ccr.jussieu.fr \\
$^{**}$Laboratoire de Math\'ematiques (UMR 6056), Universit\'e de
Reims Champagne-Ardenne and LSTA-Universit\'e Paris 6, e-mail:
amor.keziou@upmc.fr}
\date{December 22, 2007}
\maketitle

\begin{abstract}
We introduce estimation and test procedures through divergence
optimization for discrete or continuous parametric models. This
approach is based on a new dual representation for divergences. We
treat point estimation and tests for simple and composite
hypotheses, extending maximum likelihood technique. An other view
at the maximum likelihood approach, for estimation and test, is
given. We prove existence and consistency of the proposed
estimates. The limit laws of the estimates and test statistics
(including the generalized likelihood ratio one) are given both
under the null and the alternative hypotheses, and approximation
of the power functions is deduced. A new procedure of construction
of confidence regions, when the parameter may be a boundary value
of the parameter space, is proposed. Also, a solution to the
irregularity problem of the generalized likelihood ratio test
pertaining to the number of components in a mixture is given, and
a new test is proposed, based on $\chi ^{2}$-divergence on signed
finite measures and duality technique. \vspace{2mm} \newline Key
words: Parametric estimation; Parametric test; Maximum likelihood;
Mixture; Boundary valued parameter; Power function; Duality; $\phi
$-divergence.\newline
\end{abstract}


\section{Introduction and notation}

\noindent Let $\left(\mathcal{X},\mathcal{B}\right)$ be a
measurable space and $P$ be a given probability measure (p.m.) on
$\left(\mathcal{X},\mathcal{B}\right)$. Denote $\mathcal{M}$ the
real vector space of all signed finite measures on
$\left(\mathcal{X},\mathcal{B}\right)$ and $\mathcal{M}(P)$ the
vector subspace of all signed finite measures absolutely
continuous (a.c.) with respect to (w.r.t.) $P$. Denote also
$\mathcal{M}^1$ the set of all p.m.'s on
$\left(\mathcal{X},\mathcal{B}\right)$ and $\mathcal{M}^1(P)$ the
subset of all p.m.'s a.c. w.r.t. $P$. Let $\phi$ be a
proper\footnote{We say a function is proper if its domain is non
void.} closed\footnote{The closedness of $\phi$ means that if
$a_\phi$ or $b_\phi$ are finite numbers then $\phi(x)$ tends to
$\phi(a_\phi)$ or $\phi(b_\phi)$ when $x \downarrow a_\phi$ or $x
\uparrow b_\phi$, respectively.} convex function from $]-\infty,
+\infty[$ to $[0,+\infty]$ with $\phi(1)=0$ and such that its
domain $\text{dom}\phi := \left\{ x\in\mathbb{R} \text{ such that
} \phi (x)< \infty \right\}$ is an interval with endpoints $a_\phi
< 1 < b_\phi$ (which may be finite or infinite). For any signed
finite measure $Q$ in $\mathcal{M}(P)$, the $\phi$-divergence
between $Q$ and $P$ is defined by
\begin{equation}
D_\phi(Q,P):=\int_\mathcal{X} \phi\left(\frac{dQ}{dP}(x)\right)~dP(x).
\label{divRusch}
\end{equation}
When $Q$ is not a.c. w.r.t. $P$, we set $D_\phi(Q,P)=+\infty$. The
$\phi$-divergences were introduced by \cite{Csiszar1963} as
``$f$-divergences''. For all p.m. $P$, the mappings $Q\in
\mathcal{M}\mapsto D_\phi(Q,P)$ are convex and take nonnegative
values. When $Q=P$ then $D_\phi(Q,P)=0$. Furthermore, if the
function $x\mapsto\phi(x)$ is strictly convex on a neighborhood of
$x=1$, then the following fundamental property holds
\begin{equation}
D_\phi(Q,P)=0~\text{ if and only if }~Q=P.  \label{p.f.}
\end{equation}
All these properties are presented in \cite%
{Csiszar1963,Csiszar1967b,Csiszar1967} and \cite{Liese-Vajda1987}
chapter 1, for $\phi$-divergences defined on the set of all p.m.'s
$\mathcal{M}^1$. When the $\phi$-divergences are defined on
$\mathcal{M}$, then the same properties hold. Let us conclude
these few remarks quoting that in general $D_\phi(Q,P)$ and
$D_\phi(P,Q)$ are not equal. Hence, $\phi$-divergences usually are
not distances, but they merely measure some difference between two
measures. Of course a main feature of divergences between
distributions of random variables $X$ and $Y$ is the invariance
property with respect to common smooth change of variables.

\subsection{Examples of $\protect\phi$-divergences.}

\noindent When defined on $\mathcal{M}^1$, the Kullback-Leibler
$(KL)$, modified Kullback-Leibler $(KL_m)$, $\chi^{2}$, modified
$\chi^{2}$ $(\chi_{m}^{2})$, Hellinger $(H)$, and $L_{1}$
divergences are respectively associated to the convex functions
$\phi(x)=x\log x-x+1$, $\phi(x)=-\log x+x-1$,
$\phi(x)=\frac{1}{2}{(x-1)}^{2}$,
$\phi(x)=\frac{1}{2}{(x-1)}^{2}/x$, $\phi(x)=2{(\sqrt{x}-1)}^{2}$
and $\phi(x)=\left\vert x-1\right\vert$. All these divergences
except the $L_{1}$ one, belong to the class of the so called
``power divergences'' introduced in \cite{Cressie-Read1984} (see
also \cite{Liese-Vajda1987} chapter 2). They are defined through
the class of convex functions
\begin{equation}  \label{gamma convex functions}
x\in ]0,+\infty[ \mapsto\phi_{\gamma}(x):=\frac{x^{\gamma}-\gamma
x+\gamma-1}{\gamma(\gamma-1)}
\end{equation}
if $\gamma\in\mathbb{R}\setminus \left\{0,1\right\}$,
$\phi_{0}(x):=-\log x+x-1$ and $\phi_{1}(x):=x\log x-x+1$. (For
all $\gamma\in\mathbb{R}$, we define
$\phi_\gamma(0):=\lim_{x\downarrow 0}\phi_\gamma (x)$). So, the
$KL$-divergence is associated to $\phi_1$, the $KL_m$ to $\phi_0$,
the $\chi^2$ to $\phi_2$, the $\chi^2_m$ to $\phi_{-1}$ and the
Hellinger distance to $\phi_{1/2}$.\newline

\noindent We extend the definition of the power divergences functions $Q\in
\mathcal{M}^1\mapsto D_{\phi_{\gamma}}(Q,P)$ onto the whole vector space of
all signed finite measures $\mathcal{M}$ via the extension of the definition
of the convex functions $\phi_\gamma$ : For all $\gamma \in\mathbb{R}$ such
that the function $x\mapsto \phi_{\gamma }(x)$ is not defined on $]-\infty,
0[$ or defined but not convex on whole $\mathbb{R}$, set
\begin{equation}  \label{gamma convex functions sur R}
x\in ]-\infty,+\infty[\mapsto\left\{
\begin{array}{lll}
\phi_\gamma (x) & \text{ if } & x\in [0,+\infty[, \\
+\infty & \text{ if } & x\in ]-\infty,0[.
\end{array}
\right.
\end{equation}
Note that for the $\chi^2$-divergence, the corresponding $\phi$
function $\phi_2(x):=\frac{1}{2}(x-1)^2$ is defined and convex on
whole $\mathbb{R}$.
\newline

\noindent In this paper, we are interested in estimation and test
using $\phi$-divergences. An i.i.d. sample $X_1,\ldots,X_n$ with
common unknown distribution $P$ is observed and some p.m. $Q$ is
given. We aim to estimate $D_\phi(Q,P)$ and, more generally,
$\inf_{Q\in\Omega}D_\phi(Q,P)$ where $\Omega$ is some set of
measures, as well as the measure $Q^*$ achieving the infimum on
$\Omega$. In the parametric context, these problems can be well
defined and lead to new results in estimation and tests, extending
classical notions.

\subsection{Statistical examples and motivations}

\subsubsection{Tests of fit}

Let $Q_0$ and $P$ be two p.m.'s with same support $S$. Introduce a
finite partition $A_1,\ldots,A_k$ of $S$ (when $S$ is finite this
partition is the support of $Q_0$). The quantization method
consists in approximating $D_\phi(Q_0,P)$ by
$\sum_{j=1}^{k}\phi\left(\frac{Q_0(A_j)}{P(A_j)} \right)P(A_j)$
which is estimated by
\begin{equation*}
\widetilde{D_\phi}(Q_0,P)=\sum_{j=1}^{k}\phi\left(
\frac{Q_0(A_j)}{P_n(A_j)} \right)P_n(A_j),
\end{equation*}
where $P_n$ is the empirical measure associated to the data. In this vein,
goodness of fit tests have been proposed by \cite%
{ZografosFerentinosPapaionnou1990} for fixed number of classes,
and by \cite{MenendezMoralesPardoVajda1998} and
\cite{Gyorfi-Vajda2002} when the number of classes depends on the
sample size. We refer to \cite{Pardo2006} which treats these
problems extensively and contains many more references.

\subsubsection{Parametric estimation and tests}

Let $\left\{P_\theta;\theta\in\Theta\right\}$ be some parametric
model with $\Theta$ a set in $\mathbb{R}^d$. On the basis of an
i.i.d. sample $X_1,\ldots,X_n$ with distribution $P_{\theta_T}$,
we want to estimate $\theta_T$, the unknown true value of the
parameter and perform statistical tests on the parameter using
$\phi$-divergences. When all p.m.'s $P_\theta$ share the same
finite support $S$, \cite{Liese-Vajda1987}, \cite{Lindsay1994} and
\cite{Morales-Pardo-Vajda1995} introduced the so-called ``Minimum
$\phi$-divergences estimates'' (M$\phi$DE's) (Minimum Disparity
Estimators in \cite{Lindsay1994}) of the parameter $\theta_T$,
defined by
\begin{equation}  \label{MphiE}
\widetilde{\theta}_\phi:=\arg\inf_{\theta\in\Theta}D_\phi(P_{\theta},P_n).
\end{equation}
Various parametric tests can be performed based on the previous
estimates of $\phi$-divergences; see \cite{Lindsay1994} and
\cite{Morales-Pardo-Vajda1995}. The class of estimates
(\ref{MphiE}) contains the maximum likelihood estimate (MLE).
Indeed, when $\phi(x)=\phi_0(x)=-\log x+x-1$, we obtain
\begin{equation}
\widetilde{\theta}_{KL_m}:=\arg\inf_{\theta\in\Theta}KL_m(P_{\theta},P_n)
=\arg\inf_{\theta\in\Theta}\sum_{j\in S}-\log(P_\theta(j))P_n(j)=MLE.
\end{equation}
The M$\phi$DE's (\ref{MphiE}) are motivated by the fact that a
suitable choice of the divergence may lead to an estimate more
robust than the ML one (see e.g. \cite{Lindsay1994},
\cite{Basu-Lindsay1994} and \cite{JimenezShao2001}).\newline

\noindent When interested in testing hypotheses
$\mathcal{H}_{0}:\theta _{T}=\theta _{0}$ against alternatives
$\mathcal{H}_{1}:\theta _{T}\neq \theta _{0}$, where $\theta _{0}$
is a given value, we can use the statistic $D_{\phi }(P_{\theta
_{0}},P_{n})$, the plug-in estimate of the $\phi$-divergence
between $P_{\theta _{0}}$ and $P_{\theta _{T}}$, rejecting
$\mathcal{H}_{0}$ for large values of the statistic; see e.g.
\cite{Cressie-Read1984}. In the case when $\phi (x)=-\log x+x-1$,
the corresponding test based on $KL_{m}(P_{\theta _{0}},P_{n})$
does not coincide with the generalized likelihood ratio one, which
defined through the generalized likelihood ratio (GLR) $\lambda
_{n}:=2\log \frac{\sup_{\theta \in \Theta
}\prod_{i=1}^{n}p_{\theta }(X_{i})}{\prod_{i=1}^{n}p_{\theta
_{0}}(X_{i})}$. The new estimate $\widehat{KL}_{m}(P_{\theta
_{0}},P_{\theta _{T}})$ of $KL_{m}(P_{\theta _{0}},P_{\theta
_{T}})$, which is proposed in this paper, leads to the generalized
likelihood ratio test; see remark \ref{mlrt and klm div} below.
\newline

\noindent When the support $S$ is continuous, the plug-in
estimates (\ref{MphiE}) are not well defined;
\cite{Basu-Lindsay1994} investigate the so-called ``minimum
disparity estimators'' (MDE's) for continuous models, through some
common kernel smoothing method of $P_{n}$ and $P_{\theta }$. When
$\phi (x)=-\log x+x-1$, this estimate clearly, due to smoothing,
does not coincide generally with the ML one. Also, the test based
on the associated estimator of the $KL_{m}(P_{\theta
_{0}},P_{\theta _{T}})$ is different from the generalized
likelihood ratio one. Further, their estimates poses the problem
of the choice of the kernel and the window. For Hellinger
distance, see \cite{Beran1977}. For nonparametric goodness-of-fit
test, \cite{Berlinet-Vajda-vanderMeulen1998}, \cite{Berlinet1999}
proposed a test based on the estimation of the $KL_{m}$-divergence
using the smoothed kernel estimate of the density. The extension
of their results to other divergences remains an open problem; see
\cite{Berlinet1999}, \cite{Gyorfi-Liese-Vajda-Meulen1998}, and
\cite{Berlinet-Vajda-vanderMeulen1998}. All those tests are stated
for simple null hypotheses; the case of composite null hypotheses
seems difficult to handle by the above technique. In the present
paper, we treat this problem in the parametric setting. \newline

\noindent When the support $S$ is discrete infinite or continuous,
then the plug-in estimate $D_\phi(P_{\theta},P_{n})$ usually takes
infinite value when no use is done of some partition-based
approximation. In \cite{Broniatowski2002}, a new estimation
procedure is proposed in order to estimate the $KL$-divergence
between some set of p.m.'s $\Omega$ and some p.m. $P$, without
making use of any partitioning nor smoothing, but merely making
use of the well known ``dual'' representation of the
$KL$-divergence as the Fenchel-Legendre transform of the moment
generating function. Extending the paper by
\cite{Broniatowski2002}, we will use the new dual representations
of $\phi$-divergences (see \cite{B-K2004} theorem 4.4 and
\cite{Keziou2003} theorem 2.1) to define the minimum
$\phi$-divergence estimates in both discrete and continuous
parametric models. These representations are the starting point
for the definition of estimates of the parameter $\theta_T$, which
we will call ``minimum dual $\phi$-divergence estimates''
(MD$\phi$DE's). They are defined in parametric models
$\left\{P_{\theta};\theta\in\Theta\right\}$, where the p.m.'s
$P_{\theta}$ do not necessarily have finite support; it can be
discrete or continuous, bounded or not. Also the same
representations will be applied in order to estimate
$D_\phi(P_{\theta_0},P_{\theta_T})$ and
$\inf_{\theta\in\Theta_{0}}D_\phi(P_{\theta},P_{\theta_T})$ where
$\theta_0$ is a given value in $\Theta$ and $\Theta_0$ is a given
subset of $\Theta$, which leads to various simple and composite
tests pertaining to $\theta_T$, the true unknown value of the
parameter. When $\phi(x)=-\log x+x-1$, the MD$\phi$D estimate
coincides with the maximum likelihood one (see remark \ref{un
autre point de vue} below); since our approach includes also test
procedures, it will be seen that with this peculiar choice for the
function $\phi$, we recover the classical likelihood ratio test
for simple hypotheses and for composite hypotheses (see remark
\ref{mlrt and klm div} and remark \ref{mlrt and klm div c} below).
A similar approach has been proposed by \cite{LieseVajda2006}; see
their formula (118).
\newline

\noindent In any case, an exhaustive study of M$\phi$DE's seems necessary,
in a way that would include both the discrete and the continuous support
cases. This is precisely the main scope of this paper.\newline

\noindent The remainder of this paper is organized as follows. In
section 2, we recall the dual representations of
$\phi$-divergences obtained by \cite{B-K2004} theorem 4.4,
\cite{Broniatowski-Keziou2004} theorem 2.4 and \cite{Keziou2003}
theorem 2.1. Section 3 presents, through the dual representation of $\phi$%
-divergences, various estimates and tests in the parametric
framework and deals with their asymptotic properties both under
the null and the alternative hypotheses. The existence and
consistency of the proposed estimates are proved using similar
arguments as developed in \cite{Qin-Lawless1994} lemma 1. We use
the limit laws of the proposed test statistics, in a similar way
to \cite{Morales-Pardo2001}, to give an approximation to the power
functions of the tests (including the GLR one). Observe that the
power functions of the likelihood ratio type tests are not
generally known; one of our contributions is to provide explicit
power functions in the general case for simple or composite
hypotheses. As a by-product, we obtain the minimal sample size
which ensures a given power, for quite general simple or composite
hypotheses.  In section 4, we give a solution to the irregularity
problem of the GLR test of the number of components in a mixture;
we propose a new test based on the $\chi^2 $-divergence on signed
finite measures, and a new procedure of construction of confidence
regions for the parameter in the case where $\theta_T$ may be a
boundary value of the parameter space $\Theta$. All proofs are in
the Appendix. We sometimes write $Pf$ for $\int f~dP$ for any
measure $P$ and any function $f$, when defined.

\section{Fenchel Duality for $\protect\phi$-divergences}

\noindent In this section, we recall a version of the dual
representations of $\phi $-divergences obtained in \cite{B-K2004},
using Fenchel duality technique. First, we give some notations and
some results about the conjugate (or Fenchel-Legendre transform)
of real convex functions; see e.g. \cite{Rockafellar1970} for
proofs. The Fenchel-Legendre transform of $\phi $ will be denoted
$\phi ^{\ast }$, i.e.,
\begin{equation}
t\in \mathbb{R}\mapsto \phi ^{\ast }(t):=\sup_{x\in \mathbb{R}}\left\{
tx-\phi (x)\right\} ,
\end{equation}
and the endpoints of $\text{dom}\phi ^{\ast }$ (the domain of
$\phi ^{\ast }$) will be denoted $a_{\phi ^{\ast }}$ and $b_{\phi
^{\ast }}$ with $a_{\phi ^{\ast }}\leq b_{\phi ^{\ast }}$. Note
that $\phi ^{\ast }$ is a proper closed convex function. In
particular, $a_{\phi ^{\ast }}<0<b_{\phi ^{\ast }}$, $\phi
^{\ast}(0)=0$ and
\begin{equation}
a_{\phi ^{\ast }}=\lim_{y\rightarrow -\infty }\frac{\phi (y)}{y},\quad
b_{\phi ^{\ast }}=\lim_{y\rightarrow +\infty }\frac{\phi (y)}{y}.
\label{domain de varphi*}
\end{equation}
By the closedness of $\phi $, applying the duality principle, the
conjugate $\phi^{\ast \ast }$ of $\phi ^{\ast }$ coincides with
$\phi $, i.e.,
\begin{equation}
\phi ^{\ast \ast }(t):=\sup_{x\in \mathbb{R}}\left\{ tx-\phi ^{\ast
}(x)\right\} =\phi (t),~\text{ for all }t\in \mathbb{R}.
\end{equation}
For the proper convex functions defined on $\mathbb{R}$ (endowed
with the usual topology), the lower semi-continuity\footnote{We
say a function $\phi $ is lower semi-continuous if the level sets
$\left\{ x\in \mathbb{R}\text{ such that }\phi (x)\leq \alpha
\right\}$, $\alpha \in \mathbb{R}$ are all closed.} and the
closedness properties are equivalent. The function $\phi $ (resp.
$\phi ^{\ast }$) is differentiable if it is differentiable on
$]a_{\phi },b_{\phi }[$ (resp. $]a_{\phi ^{\ast }},b_{\phi^{\ast
}}[$), the interior of its domain. Also $\phi $ (resp. $\phi^{\ast
}$) is strictly convex if it is strictly convex on $]a_{\phi
},b_{\phi }[$ (resp. $]a_{\phi ^{\ast }},b_{\phi ^{\ast }}[$). The
strict convexity of $\phi $ is equivalent to the condition that
its conjugate $\phi ^{\ast }$ is essentially smooth, i.e.,
differentiable with
\begin{equation}
\begin{array}{ccccc}
\lim_{t\downarrow a_{\phi ^{\ast }}}{\phi ^{\ast }}^{\prime }(t) & = &
-\infty & \text{ if } & a_{\phi ^{\ast }}>-\infty , \\
\lim_{t\uparrow b_{\phi ^{\ast }}}{\phi ^{\ast }}^{\prime }(t) & = & +\infty
& \text{ if } & b_{\phi ^{\ast }}<+\infty .%
\end{array}
\end{equation}
Conversely, $\phi $ is essentially smooth if and only if $\phi
^{\ast }$ is strictly convex; see e.g. \cite{Rockafellar1970}
section 26 for the proofs of these properties. If $\phi $ is
differentiable, we denote $\phi ^{\prime } $ the derivative
function of $\phi $, and we define $\phi ^{\prime }(a_{\phi })$
and $\phi ^{\prime }(b_{\phi })$ to be the limits (which may be
finite or infinite) $\lim_{x\downarrow a_{\phi }}\phi ^{\prime
}(x)$ and $\lim_{x\uparrow b_{\phi }}\phi ^{\prime }(x)$,
respectively. We denote $\text{Im}\phi ^{\prime }$ the set of all
values of the function $\phi ^{\prime }$, i.e., $\text{Im}\phi
^{\prime }:=\left\{ \phi ^{\prime }(x)\text{ such that }x\in
\lbrack a_{\phi },b_{\phi }]\right\}$. If additionally the
function $\phi $ is strictly convex, then $\phi ^{\prime }$ is
increasing on $[a_{\phi },b_{\phi }]$. Hence, it is a one-to-one
function from $[a_{\phi },b_{\phi }]$ to $\text{Im}\phi ^{\prime
}$. In this case, ${\phi ^{\prime }}^{-1}$ denotes the inverse
function of $\phi ^{\prime }$ from $\text{Im}\phi ^{\prime }$ to
$[a_{\phi },b_{\phi }]$. If $\phi $ is differentiable, then for
all $x\in ]a_{\phi },b_{\phi }[$,
\begin{equation}
\phi ^{\ast }\left( \phi ^{\prime }(x)\right) =x\phi ^{\prime }(x)-\phi
\left( x\right) .  \label{forme explicite}
\end{equation}
If additionally $\phi $ is strictly convex, then for all $t\in \text{Im}\phi
^{\prime }$ we have
\begin{equation}
\phi ^{\ast }(t)=t{\phi ^{\prime }}^{-1}(t)-\phi \left( {\phi
^{\prime }}^{-1}(t)\right) \quad \text{ and }\quad {\phi ^{\ast
}}^{\prime }(t)={\phi ^{\prime }}^{-1}(t).
\end{equation}
On the other hand, if $\phi $ is essentially smooth, then the interior of
the domain of $\phi ^{\ast }$ coincides with that of $\text{Im}\phi ^{\prime
}$, i.e., $\left( a_{\phi ^{\ast }},b_{\phi ^{\ast }}\right) =\left( \phi
^{\prime }(a_{\phi }),\phi ^{\prime }(b_{\phi })\right) $.\newline

\noindent Let $\mathcal{F}$ be some class of
$\mathcal{B}$-measurable real valued functions $f$ defined on
$\mathcal{X}$, and denote $\mathcal{M}_\mathcal{F}$, the real
vector subspace of $\mathcal{M}$, defined by
\begin{equation*}
\mathcal{M}_\mathcal{F}:=\left\{Q\in \mathcal{M} \text{ such that } \int
|f|~d|Q| <\infty, \text{ for all } f\in\mathcal{F}\right\}.
\end{equation*}
In the following theorem, we recall a version of the dual
representations of $\phi$-divergences obtained by \cite{B-K2004}
(for the proof, see \cite{B-K2004} theorem 4.4).

\begin{theorem}
\label{representation duale th} Assume that $\phi$ is
differentiable. Then, for all $Q\in\mathcal{M}_\mathcal{F}$ such
that $D_\phi(Q,P)$ is finite and $
\phi^{\prime}\left(\frac{dQ}{dP}\right)$ belongs to $\mathcal{F}$,
the $\phi$-divergence $D_\phi(Q,P)$ admits the dual representation
\begin{equation}  \label{representation duale de phi}
D_\phi(Q,P)=\sup_{f\in\mathcal{F}}\left\{\int f~dQ - \int
\phi^*(f)~dP\right\},
\end{equation}
and the function $f:=\phi^{\prime}\left(\frac{dQ}{dP}\right)$ is a
dual optimal solution\footnote{i.e., the supremum in (\ref{representation duale de phi}) is achieved at
 $f:=\phi^{\prime}\left(dQ/dP\right).$}. Furthermore, if $\phi$ is essentially
smooth\footnote{Note that this is equivalent to the condition that
its conjugate $\phi^*$ is strictly convex.}, then
$f:=\phi^{\prime}\left(dQ/dP\right)$ is the unique dual optimal
solution ($P$-a.e.).
\end{theorem}

\section{Parametric estimation and tests through minimum $\protect\phi$-divergence approach and duality technique}

\noindent We consider an identifiable parametric model $\left\{
P_{\theta };\theta \in \Theta \right\} $ defined on some
measurable space $(\mathcal{X},\mathcal{B})$ and $\Theta $ is some
set in $\mathbb{R}^{d}$, not necessarily an open set. For
simplicity, we write $D_{\phi }(\theta ,\alpha )$ instead of
$D_{\phi }(P_{\theta },P_{\alpha })$. We assume that for any
$\theta $ in $\Theta $, $P_{\theta }$ has density $p_{\theta }$
with respect to some dominating $\sigma$-finite measure $\lambda$,
which can be either with countable support or not. Assume further
that the support $S$ of the p.m. $P_{\theta}$ does not depend upon
$\theta $. On the basis of an i.i.d. sample $X_{1},...,X_{n}$ with
distribution $P_{\theta _{T}}$, we intend to estimate
$\theta_{T}$, the true unknown value of the parameter, which is
assumed to be an interior point of the parameter space $\Theta$.
We will consider only strictly convex functions $\phi$ which are
essentially smooth. We will use the following assumption
\begin{equation}
\int \left| \phi ^{\prime }\left( \frac{p_{\theta }(x)}{p_{\alpha
}(x)}\right) \right| ~dP_{\theta }(x)<\infty .  \label{condition
integrabilite}
\end{equation}
Note that if the function $\phi $ satisfies
\begin{equation}
\begin{tabular}{l}
there exists $0<\delta <1$ such that for all $c$ in $\left[ 1-\delta
,1+\delta \right] $, \\
we can find numbers $c_{1},c_{2},$ $c_{3}$ such that \\
$\phi (cx)\leq c_{1}\phi (x)+c_{2}\left| x\right| +c_{3}$, for all
real $x$,
\end{tabular}
\label{condition C.0}
\end{equation}
then the assumption (\ref{condition integrabilite}) is satisfied
whenever $D_{\phi }(\theta ,\alpha )<\infty $; see e.g.
\cite{B-K2004} lemma 3.2. Also the real convex functions $\phi
_{\gamma }$ (\ref{gamma convex functions sur R}), associated to
the class of power divergences, all satisfy the condition
(\ref{condition C.0}), including all standard divergences.
\newline

\noindent For a given $\theta \in \Theta $, consider the class of functions
\begin{equation}
\mathcal{F}=\mathcal{F}_{\theta }:=\left\{ x\mapsto \phi ^{\prime
}\left( \frac{p_{\theta }(x)}{p_{\alpha }(x)}\right) ;~\alpha \in
\Theta \right\}.
\end{equation}
By application of Theorem \ref{representation duale th} above,
when assumption (\ref{condition integrabilite}) holds for any
$\alpha \in \Theta $, we obtain
\begin{equation*}
D_{\phi }(\theta ,\theta _{T})=\sup_{f\in \mathcal{F}_{\theta }}\left\{ \int
f~dP_{\theta }-\int \phi ^{\ast }(f)~dP_{\theta _{T}}\right\} ,
\end{equation*}
which, by (\ref{forme explicite}), can be written as
\begin{equation}
D_{\phi }(\theta ,\theta _{T})=\sup_{\alpha \in \Theta }\left\{
\int \phi ^{\prime }\left(\frac{p_{\theta }}{p_{\alpha }}\right)
~dP_{\theta }-\int \left[ \frac{p_{\theta }}{p_{\alpha }}\phi
^{\prime }\left( \frac{p_{\theta } }{p_{\alpha }}\right) -\phi
\left( \frac{p_{\theta }}{p_{\alpha }}\right) \right] ~dP_{\theta
_{T}}\right\}.  \label{formule de base}
\end{equation}
Furthermore, the supremum in this display is unique and it is
achieved at $\alpha =\theta _{T}$ independently upon the value of
$\theta$. Hence, it is reasonable to estimate $D_{\phi }(\theta
,\theta _{T}):=\int \phi (p_{\theta }/p_{\theta _{T}})~dP_{\theta
_{T}}$, the $\phi $-divergence between $P_{\theta }$ and
$P_{\theta _{T}}$, by
\begin{equation}
\widehat{D_{\phi }}(\theta ,\theta _{T}):=\sup_{\alpha \in \Theta
}\left\{ \int \phi ^{\prime }\left( \frac{p_{\theta }}{p_{\alpha
}}\right) ~dP_{\theta }-\int \left[ \frac{p_{\theta }}{p_{\alpha
}}\phi ^{\prime }\left( \frac{p_{\theta }}{p_{\alpha }}\right)
-\phi \left( \frac{p_{\theta } }{p_{\alpha }}\right) \right]
~dP_{n}\right\} ,  \label{estimateur de phi}
\end{equation}
in which we have replaced $P_{\theta _{T}}$ by its estimate $P_{n}$, the
empirical measure associated to the data.\newline

\noindent For a given $\theta \in \Theta $, since the supremum in
(\ref{formule de base}) is unique and it is achieved at $\alpha
=\theta _{T}$, define the following class of M-estimates of
$\theta _{T}$
\begin{equation}
\widehat{\alpha }_{\phi }(\theta ):=\arg \sup_{\alpha \in \Theta
}\left\{ \int \phi ^{\prime }\left( \frac{p_{\theta }}{p_{\alpha
}}\right) ~dP_{\theta }-\int \left[ \frac{p_{\theta }}{p_{\alpha
}}\phi ^{\prime }\left( \frac{p_{\theta }}{p_{\alpha }}\right)
-\phi \left( \frac{p_{\theta } }{p_{\alpha }}\right) \right]
~dP_{n}\right\}  \label{def teta-n(alpha)}
\end{equation}
which we call ``dual $\phi $-divergence estimates'' (D$\phi
$DE's); (in the sequel, we sometimes write $\widehat{\alpha }$
instead of $\widehat{\alpha}_{\phi }(\theta)$). Further, we have
\begin{equation*}
\inf_{\theta \in \Theta }D_{\phi }(\theta ,\theta _{T})=D_{\phi }(\theta
_{T},\theta _{T})=0.
\end{equation*}
The infimum in this display is unique and it is achieved at
$\theta =\theta_{T}$. It follows that a natural definition of
minimum $\phi$-divergence estimates of $\theta_{T}$, which we will
call ``minimum dual $\phi$-divergence estimates'' (MD$\phi $DE's),
is
\begin{equation}
\widehat{\theta }_{\phi }:=\arg \inf_{\theta \in \Theta
}\sup_{\alpha \in \Theta }\left\{ \int \phi ^{\prime }\left(
\frac{p_{\theta }}{p_{\alpha }} \right) ~dP_{\theta }-\int \left[
\frac{p_{\theta }}{p_{\alpha }}\phi ^{\prime }\left(
\frac{p_{\theta }}{p_{\alpha }}\right) -\phi \left( \frac{
p_{\theta }}{p_{\alpha }}\right) \right] ~dP_{n}\right\} .
\label{def EMphiD estimates}
\end{equation}

\noindent In order to simplify formulas (\ref{estimateur de phi}), (\ref{def
teta-n(alpha)}) and (\ref{def EMphiD estimates}), define the functions
\begin{equation}
g(\theta ,\alpha ):x\mapsto g(\theta ,\alpha ,x):=\frac{p_{\theta
}(x)}{ p_{\alpha }(x)}\phi ^{\prime }\left( \frac{p_{\theta
}(x)}{p_{\alpha }(x)} \right) -\phi \left( \frac{p_{\theta
}(x)}{p_{\alpha }(x)}\right) , \label{g}
\end{equation}
\begin{equation}
f(\theta ,\alpha ):x\mapsto f(\theta ,\alpha ,x):=\phi ^{\prime }\left(
\frac{p_{\theta }(x)}{p_{\alpha }(x)}\right)  \label{ff}
\end{equation}
and
\begin{equation}
h(\theta ,\alpha ):x\mapsto h(\theta ,\alpha ,x):=P_{\theta }f(\theta
,\alpha )-g(\theta ,\alpha ,x).  \label{h}
\end{equation}
Hence, (\ref{estimateur de phi}), (\ref{def teta-n(alpha)}) and (\ref{def
EMphiD estimates}) can be written as follows
\begin{equation}
\widehat{D_{\phi }}(\theta ,\theta _{T}):=\sup_{\alpha \in \Theta
}P_{n}h(\theta ,\alpha ),  \label{estimateur de phi simple}
\end{equation}
\begin{equation}
\widehat{\alpha }_{\phi }(\theta ):=\arg \sup_{\alpha \in \Theta
}P_{n}h(\theta ,\alpha )  \label{def teta-n(alpha) simple}
\end{equation}
and
\begin{equation}
\widehat{\theta }_{\phi }:=\arg \inf_{\theta \in \Theta }\sup_{\alpha \in
\Theta }P_{n}h(\theta ,\alpha ).  \label{def EMphiD estimates simple}
\end{equation}
Formula (\ref{formule de base}) can be written then as
\begin{equation}
D_{\phi }(\theta ,\theta _{T})=\sup_{\alpha \in \Theta }P_{\theta
_{T}}h(\theta ,\alpha ).  \label{formule de base simple}
\end{equation}
If the supremum in (\ref{def teta-n(alpha) simple}) is not unique,
we define the estimate $\widehat{\alpha }_{\phi }(\theta )$ as any
value of $\alpha \in \Theta $ that maximizes the function $\alpha
\in \theta \mapsto P_{n}h(\theta ,\alpha ).$ Also, if the infimum
in (\ref{def EMphiD estimates simple}) is not unique, the estimate
$\widehat{\theta }_{\phi }$ is defined as any value of $\theta \in
\Theta $ that minimizes the function $\theta \mapsto \sup_{\alpha
\in \Theta }P_{n}h(\theta ,\alpha ).$ Conditions assuring the
existences of the above estimates are given in  section 3.1 and
3.2 below.

\begin{remark}
For the $L_{1}$ distance, i.e. when $\phi (x)=|x-1|$, formula
(\ref{formule de base}) does not apply since the corresponding
$\phi $ function is not differentiable. However, using the general
dual representation of divergences given in \cite{B-K2004} theorem
4.1, we can obtain an explicit formula for $L_{1}$ distance
avoiding the differentiability assumption. A methodology on
estimation and testing in $L_{1}$ distance has been proposed by
\cite{DevryeLugosi2001}, and its consequences for composite
hypothesis testing and for model selection based density estimates
for nested classes of densities are presented in
\cite{DevroyeGyorfiLugosi2002} and \cite{BiauDevroye2005}.
\end{remark}

\begin{remark}
\label{un autre point de vue} (\textbf{An other view at the ML
estimate}). The maximum likelihood estimate belongs to both
classes of estimates (\ref{def teta-n(alpha) simple}) and
(\ref{def EMphiD estimates simple}). Indeed, it is obtained when
$\phi (x)=-\log x+x-1$, that is as the dual modified
$KL$-divergence estimate or as the minimum dual modified
$KL$-divergence estimate, i.e., MLE=D$\left( KL_{m}\right)
$DE=MD$\left( KL_{m}\right) $DE. Indeed, we then have $P_{\theta
}f(\theta ,\alpha )=0$ and $P_{n}h(\theta ,\alpha )=-\int \log
\left( \frac{p_{\theta }}{p_{\alpha }}\right) ~dP_{n}$. Hence by
definitions (\ref{def teta-n(alpha)}) and (\ref{def EMphiD
estimates}), we get
\begin{equation*}
\widehat{\alpha }_{KL_{m}}(\theta )=\arg \sup_{\alpha \in \Theta }-\int \log
\left( \frac{p_{\theta }}{p_{\alpha }}\right) ~dP_{n}=\arg \sup_{\alpha \in
\Theta }\int \log (p_{\alpha })~dP_{n}=MLE
\end{equation*}
independently upon $\theta$, and
\begin{equation*}
\widehat{\theta }_{KL_{m}}=\arg \inf_{\theta \in \Theta }\sup_{\alpha \in
\Theta }-\int \log \left( \frac{p_{\theta }}{p_{\alpha }}\right)
~dP_{n}=\arg \sup_{\theta \in \Theta }\int \log (p_{\theta })~dP_{n}=MLE.
\end{equation*}
So, the MLE can be seen as the estimate of $\theta _{T}$ that minimizes the
estimate of the $KL_{m}$-divergence between the parametric model $\left\{
P_{\theta };~\theta \in \Theta \right\} $ and the p.m. $P_{\theta _{T}}$.
\end{remark}

\subsection{The asymptotic properties of the D$\protect\phi$DE's
$\widehat{\protect\alpha}_\protect\protect\phi(\protect\theta)$
and
$\widehat{D_\protect\protect\phi}(\protect\theta,\protect\theta_T)$
for a given $\protect\theta$ in $\Theta$}

 \noindent This section
deals with the asymptotic properties of the estimates
(\ref{estimateur de phi simple}) and (\ref{def teta-n(alpha)
simple}). We will use similar arguments as developed in
\cite{vanderVaart1998} section 5.2 and 5.6 under classical
conditions, for the study of M-estimates. In the sequel, we assume
that condition (\ref{condition integrabilite}) holds for any
$\alpha \in \Theta $, and use $\Vert .\Vert $ to denote the
Euclidean norm in $\mathds{R}^{d}.$

\subsubsection{Consistency}
Consider the following conditions
\begin{enumerate}
 \item[(c.1)] The estimate $\widehat{\alpha}_\phi(\theta)$  exists;
 \item[(c.2)] $\sup_{\alpha\in\Theta}\left|P_{n}h(\theta,\alpha)-
P_{\theta_T}h(\theta,\alpha)\right|$ converges to zero a.s. (resp.
in probability);
 \item[(c.3)] for any positive $\epsilon$, there
exists some positive $\eta$ such that for all $\alpha\in\Theta$
satisfying $\|\alpha -\theta_T\|>\epsilon$ we have
$$P_{\theta_T}h(\theta,\alpha)
<P_{\theta_T}h(\theta,\theta_T)-\eta.$$
\end{enumerate}

\begin{remark} Condition (c.1) is fulfilled for example if the function $\alpha \in \Theta
\mapsto P_{n}h(\theta ,\alpha )$ is continuous and $\Theta $ is
compact. Condition (c.2) is satisfied if $\left\{x\mapsto
h(\theta,\alpha,x);~\alpha\in\Theta\right\}$ is a
Glivenko-Cantelli class of functions. Condition (c.3) means that
the maximizer $\alpha=\theta _T$ of the function $\alpha\mapsto
P_{\theta_T}h(\theta,\alpha)$ is well-separated. This condition
holds, for example, when the function $\alpha\in\Theta\mapsto
P_{\theta_T}h(\theta,\alpha)$ is strictly concave and $\Theta$ is
convex, which is the case for the following two examples:
\end{remark}

\begin{example}\label{example 1}
Consider the case $\phi(x)=-\log x+x-1$ and the normal model
$$\left\{\mathcal{N}(\alpha,1);~\alpha\in
\Theta=\mathbb{R}\right\}.$$ Hence, we obtain
\begin{equation}\label{eq exemple 1}
P_{\theta_T}h(\theta,\alpha)=\frac{1}{2}(\theta-\theta_T)^2-\frac{1}{2}(\alpha-\theta_T)^2.
\end{equation}
We see that condition (c.3) is satisfied; we can choose
$\eta=\frac{\epsilon^2}{2}.$
\end{example}

\begin{example}\label{example 2}
Consider the case $\phi(x)=-\log x+x-1$ and the exponential model
$$\left\{p_\alpha(x)=\alpha\exp(-\alpha x);~\alpha\in
\Theta=\mathbb{R}_+^*\right\}.$$ Hence, we obtain
\begin{equation}\label{eq exemple 2}
P_{\theta_T}h(\theta,\alpha)=-\log\theta
+\frac{\theta}{\theta_T}+\log \alpha-\frac{\alpha}{\theta_T},
 \end{equation}
which is strictly concave (in $\alpha$). Hence, condition (c.3) is
satisfied.
\end{example}

\noindent
\begin{proposition}\label{prop1}
\begin{enumerate}
\item[(1)] Under assumption (c.1-2), the estimate
$\widehat{D_\phi}(\theta,\theta_T)$ converges a.s.   (resp. in
probability) to $D_\phi(\theta,\theta_T)$.

\item[(2)] Assume that the assumptions (c.1-2-3) hold. Then the
estimate $\widehat{\alpha}_\phi(\theta)$ converges  in probability
to $\theta_T.$
\end{enumerate}
\end{proposition}


\subsubsection{Asymptotic Normality}
Assume that $\theta_T$ is an interior point of $\Theta$, the
convex function $\phi$ has continuous derivatives up to 4th order,
and the density $p_{\alpha}(x)$ has continuous partial derivatives
up to 3th order (for all $x$ $\lambda-a.e$). Denote $I_{\theta_T}$
the Fisher information matrix
\begin{equation*}
I_{\theta_T}:=\int
\frac{p_{\theta_T}^{\prime}{p_{\theta_T}^{\prime}}^T}{
p_{\theta_T}} ~ d\lambda.
\end{equation*}
In the following theorem, we  give the limit laws of the estimates
$\widehat{\alpha}_\phi(\theta)$ and $\widehat{D_\phi}
(\theta,\theta_T)$.  We will use the following assumptions.

\begin{enumerate}
 \item [(A.0)] The estimate $\widehat{\alpha}_\phi(\theta)$ exists and
 is consistent;
 \item[(A.1)] There exists a neighborhood $N(\theta_T)$
of $\theta_T$ such that the first and second order partial
derivatives (w.r.t $\alpha$) of $f(\theta,\alpha,x)p_\theta (x)$
are dominated on $N(\theta_T)$ by some $\lambda$-integrable
functions. The third order partial derivatives (w.r.t $\alpha$) of
$h(\theta,\alpha,x)$ are dominated on $N(\theta_T)$ by some
$P_{\theta_T}$-integrable functions;
 \item[(A.2)] The integrals $P_{\theta_T}
\left\|(\partial/\partial\alpha)h(\theta,\theta_T)\right\|^2$ and
$P_{\theta_T}\left\|(\partial^2/\partial\alpha^2)h(\theta,\theta_T)\right\|$
are finite, and the matrix
$P_{\theta_T}(\partial^2/\partial\alpha^2)h(\theta,\theta_T)$ is
non singular;
 \item[(A.3)] The integral
$P_{\theta_T}h(\theta,\theta_T)^2$ is finite.
\end{enumerate}

\begin{theorem} \label{th asym 1} Assume that assumptions (A.0-1-2) hold. Then, we
have
\begin{enumerate}
\item[(a)]
$\sqrt{n}\left(\widehat{\alpha}_\phi(\theta)-\theta_T\right)$
converges in distribution to a centered multivariate normal random
variable with covariance matrix
\begin{equation}  \label{variance limite 1}
V_\phi(\theta,\theta_T)=S^{-1}MS^{-1}
\end{equation}
with
$S:=-P_{\theta_T}(\partial^2/\partial\alpha^2)h(\theta,\theta_T)$
and $M:=P_{\theta_T}(\partial/\partial\alpha)h(\theta,\theta_T)
(\partial/\partial\alpha)^Th(\theta,\theta_T)$.\newline If
$\theta_T=\theta$, then
$V_\phi(\theta,\theta_T)=V(\theta_T)=I_{\theta_T}^{-1}$.

\item[(b)] If $\theta_T=\theta$, then the statistic
$\frac{2n}{\phi^{\prime\prime}(1)}\widehat{D_\phi}(\theta,\theta_T)$
converges in distribution to a $\chi^2$ random variable with $d$
degrees of freedom.

\item[(c)] If additionally assumption (A.3) holds, then when
$\theta \neq\theta_T$, we have \newline
$\sqrt{n}\left(\widehat{D_\phi}(\theta,\theta_T)-D_\phi(\theta,\theta_T)\right)$
converges in distribution to a centered normal random variable
with variance
\begin{equation}
\sigma^2_{\phi}(\theta,\theta_T)=P_{\theta_T}h(\theta,\theta_T)^2
-\left(P_{\theta_T}h(\theta,\theta_T)\right)^2.
\end{equation}
\end{enumerate}
\end{theorem}


\begin{remark}\label{motiv max loc}{\rm
\noindent Our first result (proposition \ref{prop1} above)
provides a general solution for the consistency of the global
maximum (\ref{def teta-n(alpha) simple}) under strong but usual
conditions, also difficult to be checked; see
\cite{vanderVaart1998} chapter 5. Moreover, in practice, the
optimization in (\ref{def teta-n(alpha) simple}) is handled
through gradient descent algorithms, depending on some initial
guess $\alpha_{0}\in \Theta $, which may provide a local maximum
(not necessarily global) of $P_{n}h(\theta,.)$. Hence, it is
desirable to prove that in a ``neighborhood'' of $\theta _{T}$
there exits a maximum of $P_{n}h(\theta,.)$ which indeed converges
to $\theta_{T}$; this is the scope of theorem \ref{th asym 1-2},
in the following subsection, which states that for some ``good''
$\alpha _{0}$ (near $\theta _{T}$) the algorithm provides a
consistent estimate. It is well known that, in various classical
models, the global maximizer of the likelihood function may not
exist or be inconsistent. Typical examples are provided in mixture
models. Consider the Beta-mixture model given in
\cite{Ferguson1982} section 3
$$p_\theta(x)=\theta g(x|1,1)+(1-\theta)g(x|\gamma(\theta),\beta(\theta)),$$
where $\Theta=[1/2,1]$, $g(x|\gamma(\theta),\beta(\theta))$ is the
$Be(\gamma,\beta)$-density and
$\gamma(\theta)=\theta\delta(\theta)$ and
$\beta(\theta)=(1-\theta)\delta(\theta)$ with $\delta(\theta)\to
+\infty$ sufficiently fast as $\theta\to 1.$ The ML estimate
converges to $1$ (a.s.) whatever the value of $\theta_T$ in
$\Theta;$ see \cite{Ferguson1982} section 3 for the proof.
However, if we take for example $\theta_T=3/4$, theorem \ref{th
asym 1-2} hereafter proves the existence and consistency of a
sequence of local maximizers under weak assumptions which hold for
this example. Other motivations for the results of theorem \ref{th
asym 1-2} are given in remark \ref{remark motiv 1} below.}
\end{remark}

\subsubsection{Existence, consistency and limit laws
of a sequence of local maxima}
We use similar arguments as developed in \cite{Qin-Lawless1994}
lemma 1. Assume that  $\theta_T$ is an interior point of $\Theta$,
the convex function $\phi$ has continuous derivatives up to 4th
order, and the density $p_{\alpha}(x)$ has continuous partial
derivatives up to 3th order (for all $x$ $\lambda-a.e$). In the
following theorem, we state the existence and the consistency of a
sequence of local maxima $\widetilde{\alpha}_\phi(\theta)$ and
$\widetilde{D_\phi} (\theta,\theta_T)$. We give also their limit
laws.
\begin{theorem}
\label{th asym 1-2} Assume that assumptions (A.1) and (A.2) hold.
Then, we have
\begin{enumerate}
\item[(a)] Let $B(\theta_T,n^{-1/3}):=\left\{\alpha\in\Theta;~
\|\alpha-\theta_T\|\leq n^{-1/3}\right\}$. Then, as $n \to\infty$,
with probability one, the function $\alpha\mapsto
P_nh(\theta,\alpha)$ attains its maximum value at some point
$\widetilde{\alpha}_\phi(\theta)$ in the interior of the ball $B$,
and satisfies $P_n(\partial/\partial
\alpha)h(\theta,\widetilde{\alpha}_\phi(\theta))=0$.

\item[(b)]
$\sqrt{n}\left(\widetilde{\alpha}_\phi(\theta)-\theta_T\right)$
converges in distribution to a centered multivariate normal random
variable with covariance matrix
\begin{equation}  \label{variance limite 1-2}
V_\phi(\theta,\theta_T)=S^{-1}MS^{-1}.
\end{equation}

\item[(c)] If $\theta_T=\theta$, then the statistic $\frac{2n}{
\phi^{\prime\prime}(1)}\widetilde{D_\phi}(\theta,\theta_T)$
converges in distribution to a $\chi^2$ random variable with $d$
degrees of freedom.

\item[(d)] If additionally assumption (A.3) holds, then when
$\theta \neq\theta_T$, we have \newline
$\sqrt{n}\left(\widetilde{D_\phi}(\theta,\theta_T)-D_\phi(\theta,\theta_T)
\right)$ converges in distribution to a centered normal random
variable with variance  $\sigma^2_{\phi}(\theta,\theta_T).$
\end{enumerate}
\end{theorem}

\begin{remark} \label{remark motiv 1} The results of this theorem are motivated by the
following statements
 \item [-] The estimates $\widetilde{\alpha}_\phi(\theta)$ can be calculated if the statistician
      disposes of some preknowledge of the true unknown parameter $\theta_T$.
  \item [-] The hypotheses are satisfied for a large class of
  parametric models for which the support does not depend upon
  $\theta$, such normal, log normal, exponential, Gamma, Beta,
  Weibull, ... etc; see for example \cite{vanderVaart1998}
  paragraph 5.43.
  \item [-] The maps $h(\theta,\alpha): x\mapsto
  h(\theta,\alpha,x)$ and $(\theta,\alpha)\mapsto P_{\theta_T}h(\theta,\alpha)$ are allowed to take the value $-\infty;$
  for example, take $\phi(x)=-\log x+x-1$, and consider the model
  $$\left\{P_\alpha =\alpha Cauchy(0)+(1-\alpha)\mathcal{N}(0,1);~ \alpha\in\Theta\right\},$$
  with $\Theta=[0,1]$ and $\theta_T=1/2$. Then,
  $P_{\theta_T}h(\theta,1)=-\infty$ for all $\theta\in ]0,1[$.
  \item [-] The theorem states both existence, consistency and
  asymptotic normality of the estimates.
  \item [-] The estimate $\widetilde{\alpha}_\phi(\theta)$ may
  exist and be consistent whereas $\widehat{\alpha}_\phi(\theta)$ does not in many cases.
  \item [-] One interesting situation also is if the map
 $\alpha\in\Theta\mapsto P_nh(\theta,\alpha)=0$
  is strictly concave and $\Theta$ is convex; the estimates $\widetilde{\alpha}_\phi(\theta)$
  and $\widehat{\alpha}_\phi(\theta)$ are the same.
\end{remark}

\begin{remark}
Using theorem \ref{th asym 1} part (c), the estimate
$\widehat{D_\phi}(\theta_0,\theta_T)$ can be used to perform
statistical tests (asymptotically of level $\epsilon$) of the null
hypothesis $\mathcal{H}_0:\theta_T=\theta_0$ against the
alternative $\mathcal{H}_1:\theta_T\neq \theta_0$ for a given
value $\theta_0$. Since $D_\phi(\theta_0,\theta_T)$ is nonnegative
and takes value zero only when $\theta_T=\theta_0$, the tests are
defined through the critical region
\begin{equation}  \label{CR phi}
C_\phi(\theta_0,\theta_T):=\left\{\frac{2n}{\phi^{\prime\prime}(1)}
\widehat{D_\phi}(\theta_0,\theta_T)>q_{d,\epsilon}\right\}
\end{equation}
where $q_{d,\epsilon}$ is the $(1-\epsilon)$-quantile of the
$\chi^2$ distribution with $d$ degrees of freedom. Note that these
tests are all consistent, since
$\widehat{D_\phi}(\theta_0,\theta_T)$ are $n$-consistent estimates
of $D_\phi(\theta_0,\theta_T)=0$ under $\mathcal{H}_0$, and
$\sqrt{n}$-consistent estimate of $D_\phi(\theta_0,\theta_T)>0$
under $\mathcal{H}_1 $; see part (c) and (d) in theorem \ref{th
asym 1} above. Further, the asymptotic result (d) in theorem
\ref{th asym 1} above can be used to give approximation of the
power function $\theta_T\mapsto
\beta(\theta_T):=P_{\theta_T}\left(C_\phi(\theta_0,\theta_T)\right)$.
We obtain then the following approximation
\begin{equation}  \label{approx 1}
\beta(\theta_T)\approx
1-F_\mathcal{N}\left(\frac{\sqrt{n}}{\sigma_\phi(\theta_0,\theta_T)}
\left[\frac{\phi^{\prime\prime}(1)}{2n}q_{d,\epsilon}-D_\phi(\theta_0,\theta_T)\right]\right)
\end{equation}
where $F_\mathcal{N}$ is the cumulative distribution function of a normal
random variable with mean zero and variance one. An important application of
this approximation is the approximate sample size (\ref{sample approx 1})
below that ensures a power $\beta$ for a given alternative $\theta_T\neq
\theta_0$. Let $n_0$ be the positive root of the equation
\begin{equation*}
\beta =
1-F_\mathcal{N}\left(\frac{\sqrt{n}}{\sigma_\phi(\theta_0,\theta_T)}
\left[\frac{\phi^{\prime\prime}(1)}{2n}q_{d,\epsilon}-D_\phi(\theta_0,
\theta_T)\right]\right)
\end{equation*}
i.e.,
$n_0=\frac{(a+b)-\sqrt{a(a+2b)}}{2D_\phi(\theta_0,\theta_T)^2}$
where
$a=\sigma_\phi^2(\theta_0,\theta_T)\left[F_\mathcal{N}^{-1}(1-\beta)\right]^2$
and
$b=\phi^{\prime\prime}(1)q_{d,\epsilon}D_\phi(\theta_0,\theta_T)$.
The required sample size is then
\begin{equation}  \label{sample approx 1}
n^*=[n_0]+1
\end{equation}
where $[.]$ is used here to denote ``integer part of''.
\end{remark}

\begin{remark}
\label{mlrt and klm div} (\textbf{An other view at the generalized
likelihood ratio test and approximation of the power function
through $KL_m$ -divergence}). In the particular case of the
$KL_m$-divergence, i.e., when $\phi(x)=\phi_0(x):=-\log x+x-1$, we
obtain from (\ref{CR phi}) the critical area
\begin{equation*}
C_{KL_{m}}(\theta_0,\theta_T):=\left\{2n\sup_{\alpha\in\Theta}P_{n}\log\left(
\frac{ p_{\alpha}}{p_{\theta_0}}\right) >q_{d,\epsilon}\right\}
=\left\{2\log
\frac{\sup_{\alpha\in\Theta}\prod_{i=1}^{n}p_{\alpha}(X_{i})}{
\prod_{i=1}^{n}p_{\theta_0}(X_{i})}>q_{d,\epsilon}\right\},
\end{equation*}
which is to say that the test obtained in this case is precisely the
generalized likelihood ratio one. The power approximation and the
approximate sample size guaranteeing a power $\beta$ for a given alternative
(for the GLRT) are given by (\ref{approx 1}) and (\ref{sample approx 1}),
respectively, where $\phi$ is replaced by $\phi_0$ and $D_\phi$ by $KL_m$.
\end{remark}

\subsection{The asymptotic behavior of the MD$\protect\phi$DE's}

We now explore the asymptotic properties of the estimates
$\widehat{\theta}_\phi$ and
$\widehat{\alpha}_\phi(\widehat{\theta}_\phi)$ defined in
(\ref{def EMphiD estimates simple}) and (\ref{def teta-n(alpha)
simple}). We assume that condition (\ref{condition integrabilite})
holds for any $\alpha$, $\theta\in\Theta$.

\subsubsection{Consistency}
 We state consistency under the following
assumptions

\begin{enumerate}
 \item [(c.4)] The estimates $\widehat{\theta}_\phi$ and
$\widehat{\alpha}_\phi(\widehat{\theta}_\phi)$ exist.

 \item[(c.5)] $\sup_{\left\{\alpha,\theta\in\Theta\right\}}
\left|P_{n}h(\theta,\alpha)-P_{\theta_T}h(\theta,\alpha)\right|$
tends to $0$  in probability;

\begin{enumerate}
\item[(a)] for any positive $\epsilon$, there exists some positive
$\eta$, such that for any $\alpha$ in $\Theta$ with $\left\|\alpha
-\theta_T\right\| > \epsilon$ and for all $\theta\in\Theta$, it
holds
$P_{\theta_T}h(\theta,\alpha)<P_{\theta_T}h(\theta,\theta_T)-\eta$;

\item[(b)] there exists a neighborhood of $\theta_T$, say
$N(\theta_T)$, such that for any positive $\epsilon$, there exists
some positive $\eta$ such that for all $\alpha\in N(\theta_T)$ and
all $\theta\in\Theta$ satisfying $\Vert \theta -\theta _T\Vert >
\epsilon$, it holds
$P_{\theta_T}h(\theta_T,\alpha)<P_{\theta_T}h(\theta,\alpha
)-\eta$;
\end{enumerate}

\item[(c.6)] there exists some neighborhood $N(\theta_T)$ of
$\theta_T$ and a positive function $H$ such that for all $\alpha$
in $N(\theta_T)$, $\left\|h(\theta_T,\alpha,x)\right\| \leq H(x)$
($P_{\theta_T}$-a.s.) with $P_{\theta_T}H<\infty$.
\end{enumerate}

\begin{remark}
Condition (c.5) is fulfilled if $\left\{x\mapsto
h(\theta,\alpha);~(\theta,\alpha)\in \Theta^2\right\}$ is a
Glivenko-Cantelli class of functions. Conditions (c.5.a) and
(c.5.b) mean that the saddle-point $(\theta_T,\theta_T)$, of
$(\theta,\alpha)\in \Theta\times\Theta \mapsto
P_nh(\theta,\alpha)$, is well-separated. Note that theses two
conditions are not very restrictive, they are satisfied for
example when $\Theta$ is convex and the function
$(\theta,\alpha)\in \Theta\times\Theta \mapsto
P_nh(\theta,\alpha)$ is concave in $\alpha$ (for all $\theta$) and
convex in $\theta$ (for all $\alpha$), which is the case for
example \ref{example 1} and \ref{example 2} above, both conditions
$(c.5.a)$ and $(c.5.b)$ are satisfied; we can take
$\eta=\frac{\epsilon^2}{2}.$

\end{remark}

\begin{proposition}
\label{prop2} Assume that conditions (c.4-5-6) hold. Then,

\begin{enumerate}
\item[(1)] $\sup_{\theta\in\Theta}\Vert
\widehat{\alpha}_\phi(\theta)- \theta_T\Vert$ tends to $0$  in
probability.

\item[(2)] The MD$\phi$ estimate $\widehat{\theta}_\phi$ converges
to $ \theta_T$ ~ in probability.
\end{enumerate}
\end{proposition}


\subsection{Asymptotic normality}
Assume that $\theta_T$ is an interior point of $\Theta$, the
convex function $\phi$ has continuous derivatives up to 4th order,
and the density $p_\theta(x)$ has continuous partial derivatives
up to 3th order (for all $x$ $\lambda$-a.e.). In the following
theorem we sate the asymptotic normality of the estimates
$\widehat{\theta}_\phi$ and
$\widehat{\alpha}_\phi(\widehat{\theta}_\phi).$ We will use the
following assumptions

\begin{enumerate}
 \item [(A.4)] The estimates $\widehat{\theta}_\phi$ and
$\widehat{\alpha}_\phi(\widehat{\theta}_\phi)$ exist and are
consistent;
 \item[(A.5)] There exists a neighborhood $N(\theta_T)$ of
$\theta_T$ such that the first and second order partial
derivatives (w.r.t. $\alpha$ and $\theta$) of
$f(\theta,\alpha,x)p_\theta(x)$ are dominated on
$N(\theta_T)\times N(\theta_T)$ by $\lambda$-integrable functions.
The third partial derivatives (w.r.t. $\alpha$ and $\theta$) of
$h(\theta,\alpha,x)$ are dominated on $N(\theta_T)\times
N(\theta_T)$ by some $P_{\theta_T}$-integrable functions;
\item[(A.6)] The integrals
$P_{\theta_T}\left\|(\partial/\partial\alpha)
h(\theta_T,\theta_T)\right\|^2$,
$P_{\theta_T}\left\|(\partial/\partial\theta)h(\theta_T,\theta_T)\right\|^2$,\\
$P_{\theta_T}\left\|(\partial^2/\partial\alpha^2)h(\theta_T,\theta_T)\right\|$,
$P_{\theta_T}\left\|(\partial^2/\partial\theta^2)h(\theta_T,\theta_T)
\right\| $ and
$P_{\theta_T}\left\|(\partial^2/\partial\theta\partial\alpha)h(\theta_T,\theta_T)\right\|$
are finite, and the matrix $I_{\theta_T}$ is non singular.
\end{enumerate}

\begin{theorem}
\label{th asym 2} Assume that conditions (A.4-5-6) hold. Then,
both $\sqrt{n}\left(\widehat{\theta}_\phi-\theta_T\right)$ and $
\sqrt{n}\left(\widehat{\alpha}_\phi(\widehat{\theta}_\phi)-\theta_T\right)$
converge in distribution to a centered multivariate normal random
variable with covariance matrix $V=I_{\theta_T}^{-1}$.
\end{theorem}


\subsubsection{Existence, consistency and limit laws of a sequence of local minima-maxima}
Assume that $\theta_T$ is an interior point of $\Theta$, the
convex function $\phi$ has continuous derivatives up to 4th order,
and the density $p_\theta(x)$ has continuous partial derivatives
up to 3th order (for all $x$ $\lambda$-a.e.). In the following
theorem we sate the existence and consistency of a sequence of
local minima-maxima  $\widetilde{\theta}_\phi$ and
$\widetilde{\alpha}_\phi(\widetilde{\theta}_\phi).$ We give also
their limit laws.

\begin{theorem}
\label{th asym 2-2} Assume that conditions (A.5) and (A.6) hold.
\begin{enumerate}
\item[(a)] Let $B:=\left\{\theta\in\Theta;
~\|\theta-\theta_T\|\leq n^{-1/3}\right\}$. Then, as $n\to\infty$,
with probability one, the function $(\theta,\alpha)\mapsto
P_nh(\theta,\alpha)$ attains its min-max value at some point
$\left(\widetilde{\theta}_\phi,\widetilde{\alpha}_\phi(\widetilde{\theta
}_\phi)\right)$ in the interior of $B\times B$, and satisfies
$P_n(\partial/\partial\alpha)h \left(\widetilde{\theta}_\phi,
\widehat{\alpha}_\phi(\widetilde{\theta} _\phi)\right)=0$ and\\
$P_n(\partial/\partial\theta)h\left(\widetilde{\theta} _\phi,
\widetilde{\alpha}_\phi(\widetilde{\theta}_\phi)\right)=0$.
\item[(b)] Both
$\sqrt{n}\left(\widetilde{\theta}_\phi-\theta_T\right)$ and $
\sqrt{n}\left(\widetilde{\alpha}_\phi(\widetilde{\theta}_\phi)-\theta_T\right)$
converge in distribution to a centered multivariate normal random
variable with covariance matrix $V=I_{\theta_T}^{-1}$.
\end{enumerate}
\end{theorem}

\subsection{Composite tests by minimum $\protect\phi-$divergence}

Let $\Theta_0$ be a subset of $\Theta$. We assume that there
exists an open set $B_0\subset\mathbb{R}^{d-l}$ and mappings
$r:\Theta\to \mathbb{R}^l$ and $s:B_0\to \mathbb{R}^d$ such that
the matrices $R(\theta):=\left[\frac{
\partial}{\partial\theta_i}r(\theta)\right]$ and $S(\beta):=\left[\frac{
\partial}{\partial\beta_i}s(\beta)\right]$ exist, with elements continuous,
and are of rank $l$ and $(d-l)$, respectively,
$\Theta_0=\left\{s(\beta);~ \beta\in B_0\right\}$ and
$r(\theta)=0$ for all $\theta\in\Theta_0$. Consider the composite
null hypothesis
\begin{equation}  \label{H0 et H1 dans tests composites}
\mathcal{H}_{0}~:~\theta_T\in \Theta_0~\text{ versus }~
\mathcal{H}_{1}~:~\theta_T\in \Theta\backslash\Theta_{0}.
\end{equation}
This is equivalent to
\begin{equation*}
\mathcal{H}_{0}~:~\theta_T\in s(B_{0})~\text{ versus }~
\mathcal{H}_{1}~:~\theta_T\in \Theta \backslash s(B_{0}).
\end{equation*}
Using (\ref{formule de base simple}), the $\phi$-divergence
$D_\phi(\Theta_0,\theta_T)$, between the set of distributions
$\left\{P_\theta \text{ such that }\theta\in \Theta_0\right\}$ and
the p.m. $P_{\theta_T}$, can be written as
$D_\phi(\Theta_0,\theta_T)=\inf_{\theta\in\Theta_0}
\sup_{\alpha\in\Theta}P_{\theta_T}h(\theta,\alpha)$. Hence, it can
be estimated by
\begin{equation*}
\widehat{D_\phi}(\Theta_0,\theta_T):=\inf_{\theta\in\Theta_0}\widehat{D_\phi
}(\theta,\theta_T):=\inf_{\theta\in\Theta_0}\sup_{\alpha \in \Theta
}P_{n}h(\theta,\alpha).
\end{equation*}
\noindent We use $\widehat{D_\phi}(\Theta_0,\theta_T)$ to perform
statistical test pertaining to (\ref{H0 et H1 dans tests
composites}). Since $D_\phi(\Theta_0,\theta_T):=\inf_{\theta \in
\Theta_{0}}D_\phi \left(\theta,\theta_T\right)$ is positive under
$\mathcal{H}_{1}$ and takes value $0$ only under $\mathcal{H}_{0}$
(provided that the infimum is attained on $\Theta_0$), we reject
$\mathcal{H}_0$ whenever $\widehat{D_\phi}(\Theta_0,\theta_T)$
takes large values. The following theorem provides the limit
distribution of $\widehat{D_\phi}(\Theta_0,\theta_T)$ under the
null hypothesis $\mathcal{H}_0$.

\begin{theorem}
\label{th asym 3} Let us assume that the conditions in theorem \ref{th asym
2} are satisfied. Under $\mathcal{H}_{0}$, the statistics $\frac{2n}{\phi
^{\prime \prime }(1)}\widehat{D_\phi}(\Theta_0,\theta_T)$ converge in
distribution to a $\chi ^{2}$ random variable with $~l~$ degrees of freedom.
\end{theorem}

\noindent The following theorem gives the limit laws of the test statistics $%
\frac{2n}{\phi ^{\prime \prime }(1)}\widehat{D_\phi}(\Theta_0,\theta_T)$
under the alternative hypothesis $\mathcal{H}_1:\theta_T\in\Theta\backslash%
\Theta_0$. We will use the following assumptions.

\begin{enumerate}
\item[(C.1)] The minimum of $\theta\mapsto D_\phi(\theta,\theta_T)$ on $%
\Theta_0$ is attained at some point, say $\theta^*:=s(\beta^*)$ with $%
\beta^*\in B_0$; uniqueness then follows by strict convexity of $\phi$ and
model identifiability assumption;

\item[(C.2)] There exists a neighborhood $N(\beta^*)$ of $\beta^*$
and a neighborhood $N(\theta_T)$ of $\theta_T$ such that the first
and second order partial derivatives (w.r.t. $\alpha$ and $\beta$)
of $f(s(\beta),\alpha,x)p_{s(\beta)}(x)$ are dominated on
$N(\beta^*)\times N(\theta_T)$ by $\lambda$-integrable functions.
The third partial derivatives (w.r.t. $\beta$ and $\alpha$) of
$h(s(\beta),\alpha,x)$ are dominated on $N(\beta^*)\times
N(\theta_T)$ by some $P_{\theta_T}$-integrable functions;

\item[(C.3)] The integrals $P_{\theta_T}
\left\|(\partial/\partial\alpha)h(s(\beta^*),\theta_T)\right\|^2$,
$P_{\theta_T}\left\|(\partial/\partial\beta)h(s(\beta^*),\theta_T)\right\|^2$,
\newline
$P_{\theta_T}
\left\|(\partial^2/\partial\alpha^2)h(s(\beta^*),\theta_T)\right\|$,
$P_{\theta_T}
\left\|(\partial^2/\partial\beta^2)h(s(\beta^*),\theta_T)\right\|$
and\\ $P_{\theta_T}
\left\|(\partial^2/\partial\beta\partial\alpha)h(s(\beta^*),\theta_T)\right\|$
are finite, and the matrix
\begin{equation*}
A:=\left[%
\begin{array}{cc}
A_{11} & A_{12} \\
A_{21} & A_{22}
\end{array}
\right]
\end{equation*}
is non singular, where
$A_{11}:=P_{\theta_T}(\partial^2/\partial\beta^2)h(s(
\beta^*),\theta_T)$,
$A_{22}:=P_{\theta_T}(\partial^2/\partial\alpha^2)h(s(
\beta^*),\theta_T)$ and
$A_{12}=A_{21}^T:=P_{\theta_T}(\partial^2/\partial
\beta\partial\alpha)h(s(\beta^*),\theta_T)$.

\item[(C.4)] The integral
$P_{\theta_T}\left\|h(s(\beta^*),\theta_T)\right\|^2$ is finite.
\end{enumerate}

\noindent Denote $\widehat{\beta}_\phi$ and
$\widehat{\alpha}_\phi(\widehat{\beta}_\phi)$ the min-max optimal
solution of $$\widehat{D_\phi} (\Theta_0,\theta_T):=\inf_{\beta\in
B_0}\sup_{\alpha\in\Theta}P_nh(s(\beta),\alpha),$$ and let $
B(\beta^*,n^{-1/3}):=\left\{\beta\in B_0; \|\beta-\beta^*\|\leq
n^{-1/3}\right\}$,
$c_n:=(\widehat{\beta}_\phi^T,\widehat{\alpha}_\phi(\widehat{\beta}_\phi)^T)^T$,
$c^*:=({\beta^*}^T,\theta_T^T)^T$ and $F$ the matrix defined by
\begin{equation*}
F:=P_{\theta_T}\left[
\begin{array}{c}
(\partial/\partial\beta)h(s(\beta^*),\theta_T) \\
(\partial/\partial\alpha)h(s(\beta^*),\theta_T)
\end{array}
\right] \left[
\begin{array}{c}
(\partial/\partial\beta)h(s(\beta^*),\theta_T) \\
(\partial/\partial\alpha)h(s(\beta^*),\theta_T)
\end{array}%
\right]^T.
\end{equation*}

\begin{enumerate}
\item [(C.5)] The estimates $\widehat{\beta}_\phi$ and
$\widehat{\alpha}_\phi(\widehat{\beta}_\phi)$ exist and are
consistent estimators for $\beta^*$ and $\theta_T$ respectively.
\end{enumerate}

\begin{theorem}
\label{th asym 4} Assume that conditions (C.1-2-3-4-5) hold. Then,
under the alternative hypothesis $\mathcal{H}_1$, we have
\begin{enumerate}
\item[(a)] $\sqrt{n}\left(c_n-c^*\right)$ converges in
distribution to a centered multivariate normal random variable
with covariance matrix $V=A^{-1}FA^{-1}$.

\item[(b)] If additionally the condition (C.6) holds, then
$\sqrt{n}\left(\widehat{D_\phi}(\Theta_0,\theta_T)-D_\phi(\Theta_0,\theta_T)\right)$
converges in distribution to a centered normal random variable
with variance
\begin{equation}  \label{sigma2 c}
\sigma_\phi^2(\beta^*,\theta_T)= P_{\theta_T}h(s(\beta^*),\theta_T)^2-
\left(P_{\theta_T}h(s(\beta^*),\theta_T)\right)^2.
\end{equation}
\end{enumerate}
\end{theorem}

\begin{remark}
Using theorem \ref{th asym 3}, the estimate
$\widehat{D_\phi}(\Theta_0,\theta_T)$ can be used to perform
statistical tests (asymptotically of level $\epsilon$) of the null
hypothesis $\mathcal{H}_0:\theta_T\in\Theta_0$ against the
alternative $\mathcal{H}_1:\theta_T\in \Theta\backslash\Theta_0$.
Since $D_\phi(\Theta_0,\theta_T)$ is nonnegative and takes value
zero only when $\theta_T\in\Theta_0$, the tests are defined
through the critical region
\begin{equation}  \label{CR phi c}
C_\phi(\Theta_0,\theta_T):=\left\{\frac{2n}{\phi^{\prime\prime}(1)}
\widehat{D_\phi}(\Theta_0,\theta_T)>q_{l,\epsilon}\right\},
\end{equation}
where $q_{l,\epsilon}$ is the $(1-\epsilon)$-quantile of the
$\chi^2$ distribution with $l$ degrees of freedom. Note that these
tests are all consistent, since
$\widehat{D_\phi}(\Theta_0,\theta_T)$ are $n$-consistent estimates
of $D_\phi(\Theta_0,\theta_T)=0$ under $\mathcal{H}_0$, and
$\sqrt{ n}$-consistent estimate of $D_\phi(\Theta_0,\theta_T)>0$
under $\mathcal{H}_1 $; see theorem \ref{th asym 3} and theorem
\ref{th asym 4} part (c). Further, the asymptotic result (c) in
theorem \ref{th asym 4} above can be used to give an approximation
to the power function $\theta_T\mapsto
\beta(\theta_T):=P_{\theta_T}\left(C_\phi(\Theta_0,\theta_T)\right)$.
We obtain then the following approximation
\begin{equation}  \label{approx 1 c}
\beta(\theta_T)\approx
1-F_\mathcal{N}\left(\frac{\sqrt{n}}{\sigma_\phi(\beta^*,\theta_T)}
\left[\frac{\phi^{\prime\prime}(1)}{2n}
q_{l,\epsilon}-D_\phi(\Theta_0,\theta_T)\right]\right)
\end{equation}
where $F_\mathcal{N}$ is the cumulative distribution function of a
normal variable with mean zero and variance one. An important
application of this approximation is the approximate sample size
(\ref{sample approx 1 c}) below that ensures a power $\beta$ for a
given alternative $\theta_T\in\Theta\backslash\Theta_0$. Let $n_0$
be the positive root of the equation
\begin{equation*}
\beta =
1-F_\mathcal{N}\left(\frac{\sqrt{n}}{\sigma_\phi(\beta^*,\theta_T)}
\left[\frac{\phi^{\prime\prime}(1)}{2n}q_{l,\epsilon}-D_\phi(\Theta_0,
\theta_T)\right]\right)
\end{equation*}
i.e.,
$n_0=\frac{(a+b)-\sqrt{a(a+2b)}}{2D_\phi(\Theta_0,\theta_T)^2}$
where
$a=\sigma_\phi^2(\beta^*,\theta_T)\left[F_\mathcal{N}^{-1}(1-\beta)\right]^2$
and
$b=\phi^{\prime\prime}(1)q_{l,\epsilon}D_\phi(\Theta_0,\theta_T)$.
The required sample size is then
\begin{equation}  \label{sample approx 1 c}
n^*=[n_0]+1
\end{equation}
where $[.]$ is used here to denote ``integer part of''.
\end{remark}

\begin{remark}
\label{mlrt and klm div c} (\textbf{An other view at the
generalized likelihood ratio test for composite hypotheses, and
approximation of the power function through $KL_m$-divergence}).
In the particular case of the $KL_m$-divergence, i.e., when
$\phi(x)=\phi_0(x):=-\log x+x-1$, we obtain from (\ref{CR phi c})
the critical area
\begin{equation*}
C_{KL_{m}}(\Theta_0,\theta_T)=\left\{2\log
\frac{\sup_{\alpha\in\Theta} \prod_{i=1}^{n}p_{\alpha}(X_{i})}{
\sup_{\theta\in\Theta_0}
\prod_{i=1}^{n}p_{\theta}(X_{i})}>q_{l,\epsilon}\right\},
\end{equation*}
which is to say that the test obtained in this case is precisely the
generalized likelihood ratio test associated to (\ref{H0 et H1 dans tests
composites}). The power approximation and the approximate sample size
guaranteeing a power $\beta$ for a given alternative (for the GLRT) are
given by (\ref{approx 1 c}) and (\ref{sample approx 1 c}), respectively,
where $\phi$ is replaced by $\phi_0$ and $D_\phi$ by $KL_m$.
\end{remark}

\section{Non regular models. A simple solution for the case of mixture models}

\noindent The test problem for the number of components of a
finite mixture has been extensively treated when the total number
of components $k$ is equal to $2$, leading to a satisfactory
solution; the limit distribution of the generalized likelihood
ratio statistic is non standard, since it is $0.5\delta
_{0}+0.5\chi ^{2}(1)$, a mixture of a Dirac mass at $0$ and a
$\chi ^{2}(1)$ with weights equal to $1/2$; see e.g.
\cite{TitteringtonSmithMakov1985} and \cite{Self_Liang1987}. When
$k>2$, the problem is much more involved. \cite{Self_Liang1987}
obtained the limit distribution of the generalized likelihood
ratio statistic, which is non standard and complex. This result
yields formidable numerical difficulties for the calculation of
the critical value of the test. In section 5.1 below, we propose a
unified treatment for all these cases, with simple and standard
limit distribution both when the parameter $\theta _{T}$ is an
interior or a boundary point of the parameter space $\Theta $. On
the other hand, confidence regions for the mixture parameter
$\theta _{T}$ even when $k=2$ are intractable through the
generalized likelihood ratio statistic. Indeed, the limit law of
the generalized likelihood ratio statistic depends heavily on the
fact that $\theta $ is a boundary or an interior point of the
parameter space. For example, when $k=2$, the limit distribution
of the generalized likelihood ratio statistic is $0.5\delta
_{0}+0.5\chi ^{2}(1)$ when $\theta =0$ and $\chi ^{2}(1)$ when
$0<\theta <1$. Therefore, the confidence level is not defined
uniquely. At the opposite, we will prove in section 5.3 that the
proposed dual $\chi ^{2}$-statistic yields quite standard
confidence regions even when $k>2$.

\subsection{Notations}

Let $\left\{ P_{a_{1}}^{(1)};a_{1}\in A_{1}\right\} $, $\ldots $,
$\left\{ P_{a_{k}}^{(k)};a_{k}\in A_{k}\right\} $ be
$k$-parametric models where $A_{1},\ldots ,A_{k}$ are $k$ ($k\geq
2$)  sets in $\mathbb{R} ^{d_{1}},\ldots ,\mathbb{R}^{d_{k}}$ and
$d_{1},\ldots ,d_{k}\in \mathbb{N}^{\ast }$. Denote $P_{\theta }$
the mixture model
\begin{equation}
P_{\theta }:=\sum_{i=1}^{k}w_{i}P_{a_{i}}^{(i)}  \label{mixture model}
\end{equation}%
where $0\leq w_{i}\leq 1$, $\sum w_{i}=1$ and
\begin{equation}
\theta \in \Theta :=\left\{ (w_{1},\ldots ,w_{k},a_{1},\ldots
,a_{k})^{T}\in \lbrack 0,1]^{k}\times A_{1}\times \cdots \times
A_{k}\text{ such that } \sum_{i=1}^{k}w_{i}=1\right\} ,
\end{equation}%
and assume that the model is identifiable. Let $k_{0}\in \left\{
1,\ldots ,k-1\right\} $. We test if $(k-k_{0})$ components in
(\ref{mixture model}) have null coefficients. We assume that their
labels are $k_{0}+1,...,k.$  Denote $\Theta _{0}$ the subset of
$\Theta $ defined by
\begin{equation*}
\Theta _{0}:=\left\{ \theta \in \Theta \text{ such that
}w_{k_{0}+1}=\cdots =w_{k}=0\right\}.
\end{equation*}
On the basis of an i.i.d sample $X_{1},\ldots ,X_{n}$ with
distribution $P_{\theta _{T}}$, $\theta _{T}\in \Theta $, we
intend to perform tests of the hypothesis
\begin{equation}
\mathcal{H}_{0}:\theta _{T}\in \Theta _{0}~\text{ against the
alternative }~ \mathcal{H}_{1}:\theta _{T}\in \Theta \setminus
\Theta _{0}. \label{number components test}
\end{equation}%
It is known that the generalized likelihood ratio test, based on the
statistic
\begin{equation}
2\log \lambda :=2\log \frac{\sup_{\theta \in \Theta
}\prod_{i=1}^{n}p_{\theta }(X_{i})}{\sup_{\theta \in \Theta
_{0}}\prod_{i=1}^{n}p_{\theta }(X_{i})},  \label{GLRS_melange}
\end{equation}%
is not valid for this problem, since the asymptotic approximation
by $\chi ^{2}$ distribution does not hold in this case; the
problem is due to the fact that the null value of $\theta _{T}$ is
not in the interior of the parameter space $\Theta $. We clarify
now this problem. For simplicity, consider a mixture of two known
densities $p_{0}$ and $p_{1}$ with $p_{0}\neq p_{1}$:
\begin{equation}
p_{\theta }=(1-\theta )p_{0}+\theta p_{1}\text{ where }\theta \in \Theta
:=[0,1].  \label{mixture}
\end{equation}
Given data $X_{1},\ldots ,X_{n}$ with distribution $P_{\theta
_{T}}$, $\theta _{T}\in \lbrack 0,1]$, consider the test problem
\begin{equation}
\mathcal{H}_{0}:\theta _{T}=0~\text{ against the alternative
}~\mathcal{H}_{1}:\theta _{T}>0.  \label{simple number components
test}
\end{equation}
The generalized likelihood ratio statistic for this test problem is
\begin{equation}
W_{n}(0):=2log\frac{L(\widehat{\theta })}{L(0)},  \label{W_n}
\end{equation}
where $L(\theta ):=\prod_{i=1}^{n}\left[ (1-\theta
)p_{0}(X_{i})+\theta p_{1}(X_{i})\right] $ for all $\theta \in
\lbrack 0,1]$, and $\widehat{\theta }$ is the MLE of $\theta$.
Using the strict concavity of the function $\theta \in \lbrack
0,1]\mapsto l(\theta ):=\log L(\theta )$, it is clear that
$\widehat{\theta }=0$ whenever $l_{+}^{\prime }(0)$, the
derivative on the right at $\theta =0$ of $\theta \mapsto l(\theta
)$, is nonpositive. Hence, we can write
\begin{eqnarray}
P_{0}\left\{ W_{n}=0\right\}  &\geq &P_{0}\left\{ \widehat{\theta}
=0\right\} ~=~P_{0}\left\{ l_{+}^{\prime }(0)\leq 0\right\}
=P_{0}\left\{
\sum_{i=1}^{n}\frac{p_{0}(X_{i})}{p_{1}(X_{i})}-n\leq 0\right\}   \notag \\
&=&P_{0}\left\{ \sqrt{n}\left(
\frac{1}{n}\sum_{i=1}^{n}\frac{p_{0}(X_{i})}{
p_{1}(X_{i})}-1\right) \leq 0\right\}
\end{eqnarray}
which, by the CLT, tends to $1/2$ (if $1\neq E(Y_{i}^{2})<\infty $
where $Y_{i}:=p_{0}(X_{i})/p_{1}(X_{i})$) since the random
variables $Y_{i}$ are i.i.d with $E(Y_{i})=1$ under
$\mathcal{H}_{0}$. This proves that the convergence in
distribution of the generalized likelihood ratio statistic
$W_{n}(0)$ to a $\chi ^{2}$ random variable (under
$\mathcal{H}_{0}$) does not hold. Under suitable regularity
conditions we can prove that the limit distribution of the
statistic $W_{n}$ in (\ref{W_n}) is $0.5\delta _{0}+0.5\chi
_{1}^{2}$, a mixture of the $\chi ^{2}$-distribution and the Dirac
measure at zero; see \cite{Self_Liang1987}.

\begin{figure}[h]
\centerline{
  \begin{tabular}{ c  c }
 \includegraphics[width=.4\textwidth]{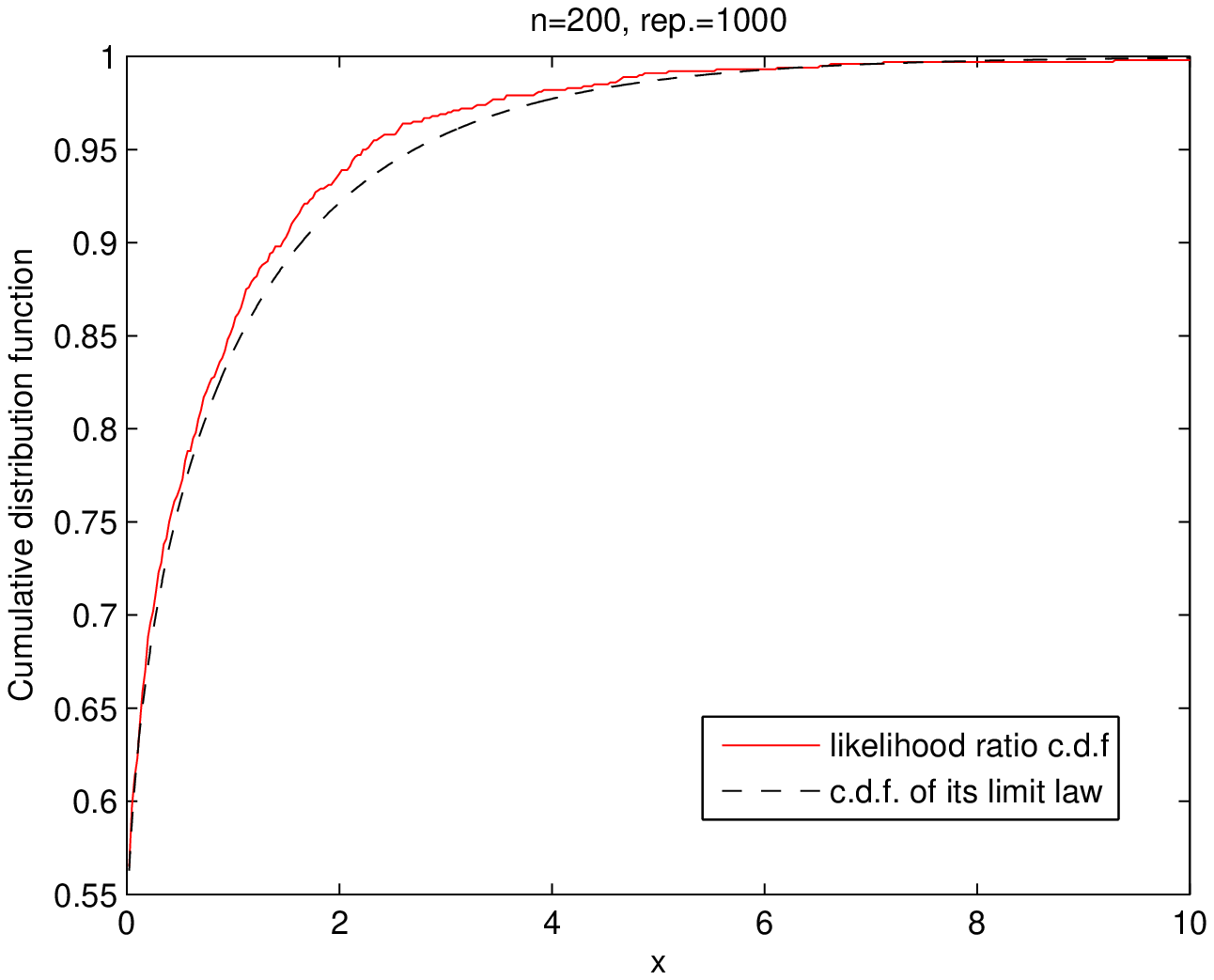}
    &   \includegraphics[width=.4\textwidth]{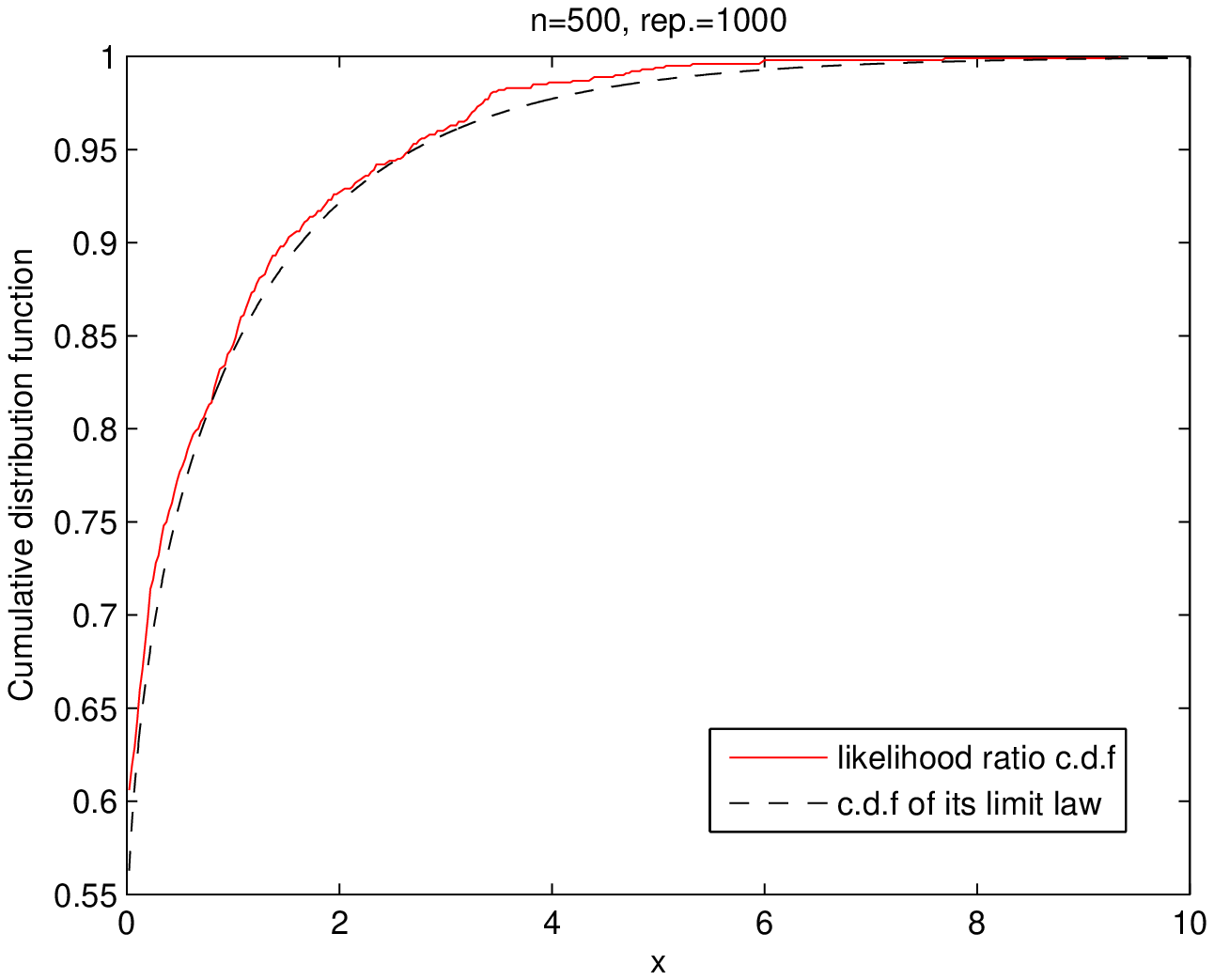}
   \end{tabular} }
\par
\begin{center}
\includegraphics[width=.4\textwidth]{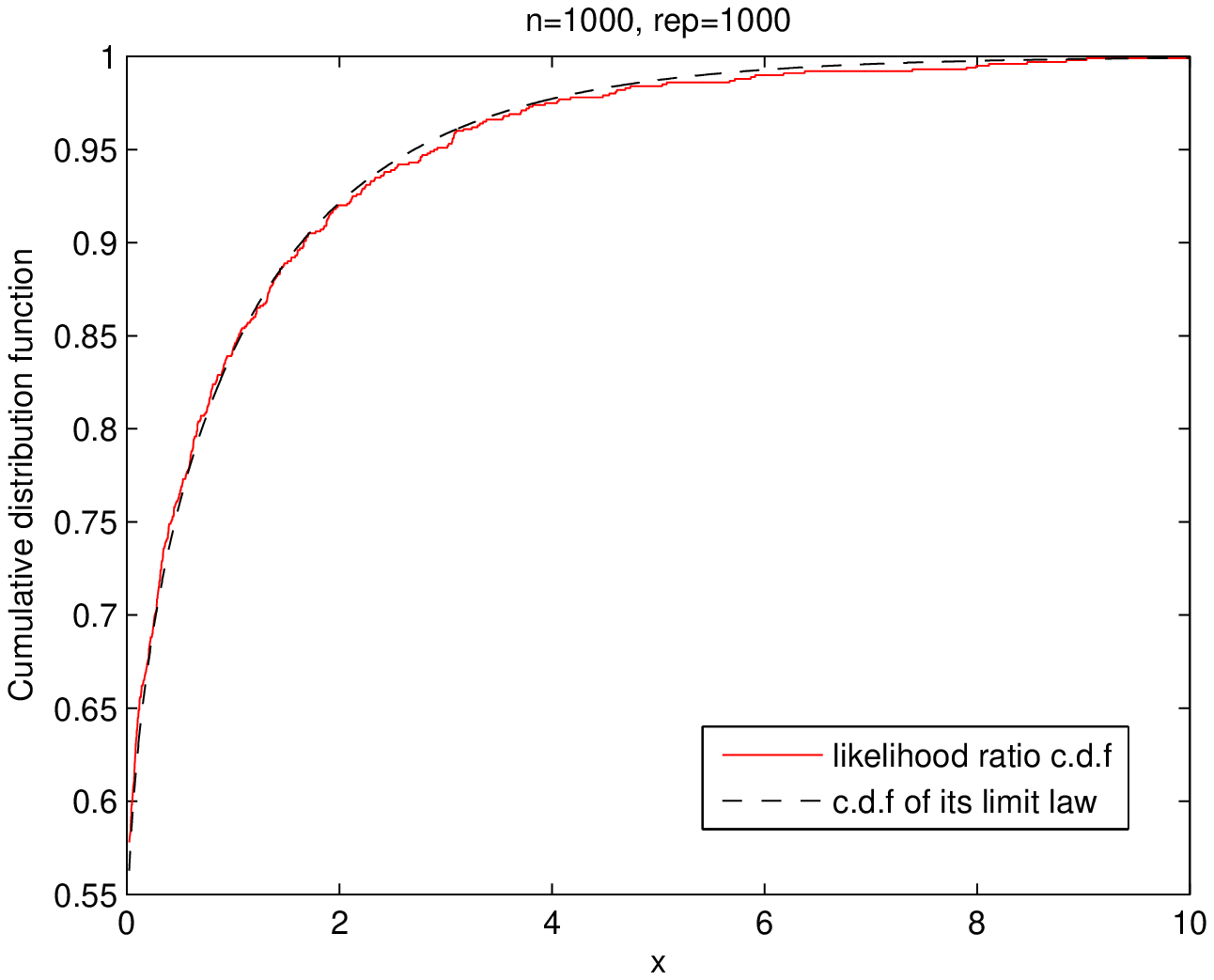}
\end{center}
\caption{Empirical distribution of the GLR and its limit
distribution} \label{the limit of the GLRS}
\end{figure}

\noindent Moreover, in the case of more than two components and
$k-k_0\geq 2$, the limit distribution of the GLR statistic
(\ref{GLRS_melange}) under $\mathcal{H}_0$ is complicate and not
standard (not a $\chi^2$ distribution) which poses some difficulty
in determining the critical value that will give correct
asymptotic size; see \cite{Self_Liang1987}. On the other hand, the
likelihood ratio statistic
\begin{equation}  \label{W_n theta}
W_n(\theta) := 2log\frac{L(\widehat{\theta})}{L(\theta)}
\end{equation}
can not be used to construct asymptotic confidence region for the
parameter $\theta_T$ since its limit law is not the same when
$\theta_T=0$ and $\theta_T>0$.

\noindent In figure 1, we illustrate the accuracy of the
approximation of distribution of the GLR by its limit $0.5\delta
_{0}+0.5\chi _{1}^{2}$; we plot the cumulative distribution
function (c.d.f) of both the limit law, and the observed GLR's
obtained from 1000 independent runs of samples with sizes $n=200$,
$n=500$ and $n=1000$, with $P_{0}=\mathcal{N}(0,1)$ and $P_{1}=
\mathcal{N}(0.5,1)$.

\subsection{A simple solution to the problem of testing the number of
components in a mixture}

We propose the following simple solution : Consider the following set of
signed finite measures
\begin{equation}
p_{\theta }:=(1-\theta )p_{0}+\theta p_{1}\text{ where }\theta \in
\mathds{R}.  \label{mixture signed measures}
\end{equation}%
This set (of signed finite measures with mass one) obviously contains the
mixture model (\ref{mixture}). In particular, the null value of $\theta _{T}$
(i.e., $\theta _{T}=0$) is an interior point of the parameter space $%
\mathds{R}$. The likelihood ratio test (for a model of signed
measures) cannot be used since the log-likelihood $l(\theta )$ may
be infinite (when $\theta <0$ or $\theta >1$). In the context of
divergences, this means that the estimate
$\widehat{KL_{m}}(P_{0},P_{\theta _{T}})$ may be infinite if we
consider the model (\ref{mixture signed measures}), which is due
to the fact that the corresponding convex function $\phi (x)=-\log
x+x-1$ is infinite on $\mathbb{R}_{-}$. This suggests to use a
divergence associated to a convex function $\phi $ which is finite
on all $\mathbb{R}$, for instance, the $\chi ^{2}$-divergence
(which is associated to the convex function $\phi (x)=
\frac{1}{2}(x-1)^{2}$). So, in order to perform a test
asymptotically of level $\epsilon $ for (\ref{simple number
components test}), we propose to use the following estimate of the
$\chi ^{2}$-divergence between $P_{0}$ and $P_{\theta _{T}}$
\begin{equation}
\widetilde{\chi ^{2}}(0,\theta _{T})=\sup_{\alpha \in \Theta
_{e}}\left\{ P_{0}f(0,\alpha )-P_{n}g(0,\alpha)\right\},
\label{estimateur de chi2 mesures signees}
\end{equation}
where $f(0,\alpha )=p_{0}/p_{\alpha }-1$ and $g(0,\alpha
)=1/2(p_{0}/p_{\alpha }+1)(p_{0}/p_{\alpha }-1)$ as a consequence of
definitions (\ref{ff}) and (\ref{g}), and $\Theta _{e}$ is the new parameter
space which we define as follows
\begin{equation*}
\Theta _{e}:=\left\{ \alpha \in \mathds{R}\text{ such that }\int \left|
f(0,\alpha )\right| ~dP_{0}\text{ is finite }\right\} .
\end{equation*}
The value of the parameter $\theta _{T}$ under the null hypothesis
$\mathcal{ H}_{0}$, i.e., $\theta _{T}=0$, is in the interior of
the new parameter space $\Theta _{e}$ which is generally non void.
Hence, under conditions of theorem \ref{th asym 1} where $\Theta $
is replaced by $\Theta _{e}$ and $\theta$ by zero, under
$\mathcal{H}_{0}$ the statistic $2n\widetilde{\chi
^{2}}(0,\theta_{T})$ converges in distribution to a $\chi ^{2}$
random variable with one degree of freedom; the critical region
takes then the form
\begin{equation}
CR:=\left\{ 2n\widetilde{\chi ^{2}}(0,\theta _{T})>q_{1,\epsilon
}\right\},
\end{equation}
where $q_{1,\epsilon }$ is the $(1-\epsilon )$-quantile of the $\chi ^{2}$
distribution with one degree of freedom. Obviously other divergences which
are associated to convex functions finite on all $\mathbb{R}$ can be used.
The use of the $\chi ^{2}$-divergence is recommended. Indeed, for regular
cases (for example for multinomial goodness-of-fit tests) $\chi ^{2}$-test
is equivalent (in Pitman sense) to the generalized likelihood ratio one; see
also \cite{Cressie-Read1984} sections 3.1 and 3.2 for other motivations in
favor of the $\chi ^{2}$ approach.

\begin{figure}[h]
\centerline{
  \begin{tabular}{ c  c }
 \includegraphics[width=.4\textwidth]{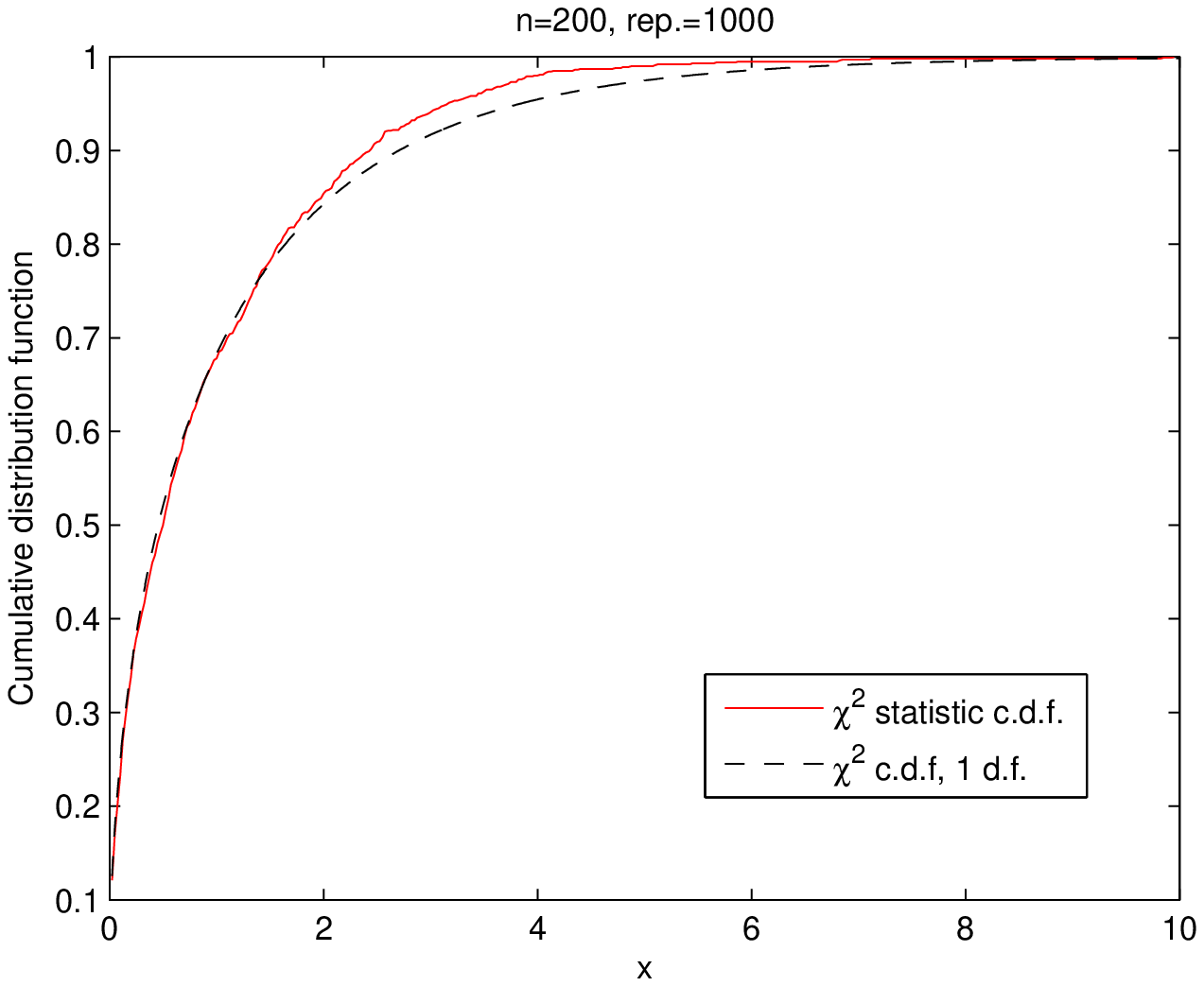}
    &   \includegraphics[width=.4\textwidth]{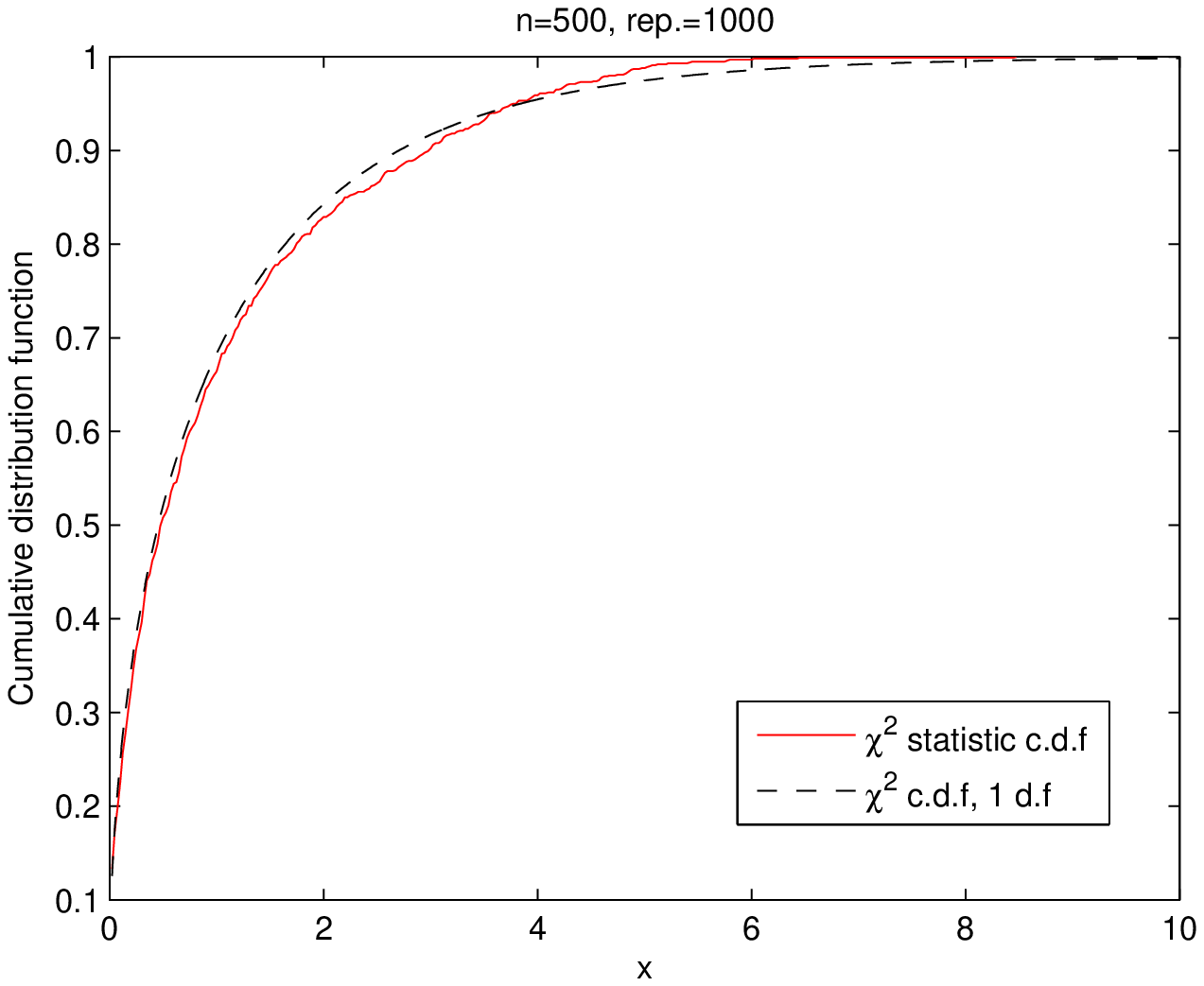}
   \end{tabular} }
\par
\begin{center}
\includegraphics[width=.4\textwidth]{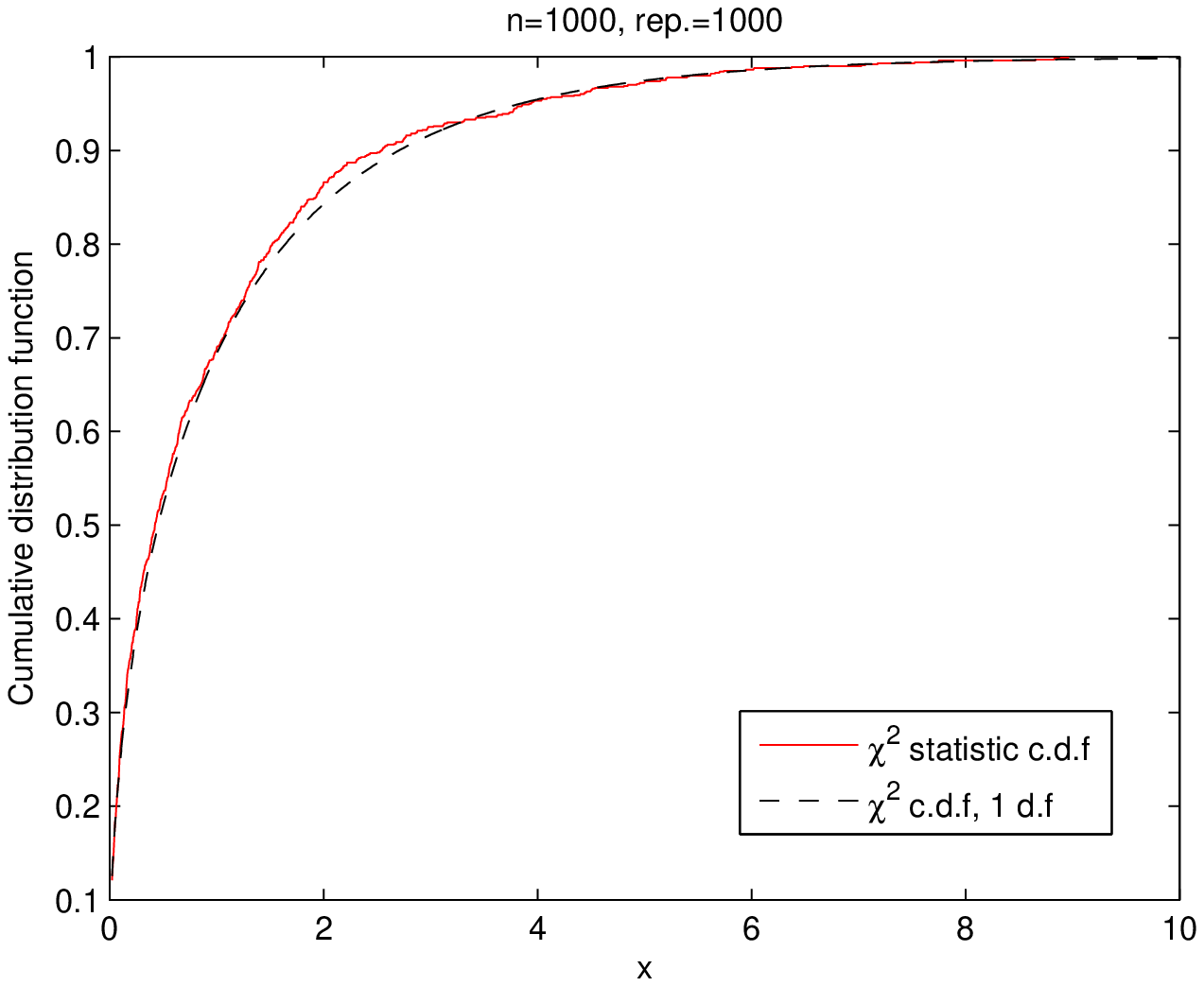}
\end{center}
\caption{Empirical distribution of the dual
$\protect\chi^{2}$-statistic and its limit law} \label{the limit
of the GLRS}
\end{figure}

\noindent In figure 2, we illustrate the accuracy of the
approximation of the distribution of the proposed dual
$\chi^{2}$-statistic by the $\chi^{2}(1)$; we plot the cumulative
distribution function (c.d.f) of both the limit law, and the dual
$\chi^{2}$-statistic obtained from 1000 independent runs of
samples with sizes $n=200$, $n=500$ and $n=1000$, with $P_{0}=
\mathcal{N}(0,1)$ and $P_{1}=\mathcal{N}(0.5,1)$. We observe that
the approximation is as satisfactory as it is in figure 1 for the
GLR  case, so that the extension of the model to signed finite
measures does not affect the quality of the approximation of the
limit distribution.

\subsection{Confidence regions for the mixture parameters}

We propose the following solution to the confidence region problem when the
parameter may be a boundary value of the parameter space: The estimate
\begin{equation}
\widetilde{\chi ^{2}}(\theta ,\theta _{T})=\sup_{\alpha \in \Theta
_{e}(\theta )}\left\{ P_{\theta }f(\theta ,\alpha )-P_{n}g(\theta ,\alpha
)\right\} ,  \label{estimateur de chi2 mesures signees theta}
\end{equation}%
where
\begin{equation*}
\Theta _{e}(\theta ):=\left\{ \alpha \in \mathds{R}\text{ such that }\int
\left| f(\theta ,\alpha )\right| ~dP_{\theta }\text{ is finite }\right\} ,
\end{equation*}%
can be used to construct asymptotic confidence region for the
parameter $\theta _{T}$ with level $(1-\epsilon )$ defined by
\begin{equation*}
C:=\left\{ \theta \in \Theta \text{ such that
}2n\widetilde{\chi^{2}} (\theta ,\theta _{T})\leq q_{1,\epsilon
}\right\}.
\end{equation*}
In fact, $\lim_{n\rightarrow \infty }P_{\theta _{T}}\left( \theta _{T}\in
C\right) =1-\epsilon $ both when $\theta _{T}=0$ or $\theta _{T}>0$ since
the statistic $2n\widetilde{\chi ^{2}}(\theta _{T},\theta _{T})$ converges
in distribution to $\chi ^{2}$ random variable with one degree of freedom
both when $\theta _{T}=0$ or $\theta _{T}>0$. We give now the form of the
critical region and the confidence region in the multivariate case, i.e., in
the case of the general model (\ref{mixture model}). For all $\theta \in
\Theta $, define the set
\begin{equation*}
\Theta _{e}(\theta ):=\left\{ \alpha \in \mathds{R}^{k}\times
A_{1}\times \cdots \times A_{k}\text{ such that
}\sum_{i=1}^{k}\alpha _{i}=1\text{ and } \int \left| f(\theta
,\alpha )\right| ~dP_{\theta }\text{ is finite } \right\} ,
\end{equation*}
and the statistic
\begin{equation*}
\widetilde{\chi ^{2}}(\Theta _{0},\theta _{T}):=\inf_{\theta \in
\Theta_{0}} \widetilde{\chi ^{2}}(\theta ,\theta
_{T}):=\inf_{\theta \in \Theta _{0}}\sup_{\alpha \in \Theta
_{e}(\theta )}\left\{ P_{\theta }f(\theta ,\alpha )-P_{n}g(\theta
,\alpha )\right\}.
\end{equation*}
Under some conditions similar to that in theorems 3.1, 3.2 and 3.3, we can
prove, under the null hypothesis $\mathcal{H}_{0}$ in (\ref{number
components test}), that the statistic $2n\widetilde{\chi ^{2}}(\Theta
_{0},\theta _{T})$ converges in distribution to $\chi ^{2}$ random variable
with $(k-k_{0})$ degrees of freedom. Also, the statistic $2n\widetilde{\chi
^{2}}(\theta ,\theta _{T})$ when $\theta =\theta _{T}$ converges in
distribution to $\chi ^{2}$ random variable with $d:=k-1+d_{1}+\cdots +d_{k}$
degrees of freedom in both case when $\theta _{T}$ is a boundary value or
not. Hence, the critical region is given by
\begin{equation*}
CR:=\left\{ 2n\widetilde{\chi ^{2}}(\Theta _{0},\theta
_{T})>q_{k-k_{0},\epsilon }\right\} ,
\end{equation*}
and
\begin{equation*}
C:=\left\{ \theta \in \Theta \text{ such that }2n\widetilde{\chi
^{2}} (\theta ,\theta _{T})\leq q_{d,\epsilon }\right\}
\end{equation*}%
is an asymptotic confidence region for $\theta _{T}$ of level
$\epsilon$ both when $\theta _{T}$ is a boundary value or not.

\subsection{Approximation of the power function of the likelihood ratio
statistic: simulation results}

In the context of the exponential model $p_{\theta }(x)=\theta
\exp\left\{\theta x\right\}$, we consider the problem of testing
\begin{equation*}
\mathcal{H}_{0}:\theta _{T}=1\quad \text{versus}\quad
\mathcal{H}_{1}:\theta_{T}\neq 1
\end{equation*}
using the GLR. We recall that the power function of the GLR test
is
\begin{equation}
\theta _{T}\mapsto \beta(\theta_T):=P_{\theta _{T}}\left\{
2n\widehat{KL_{m}}\left( 1,\theta _{T}\right) \geq
q_{1,0.5}\right\}   \label{Power}
\end{equation}
and its approximation is
\begin{equation}
\widehat{\beta}(\theta _{T})=1-F_{\mathcal{N}}\left(
\frac{\sqrt{n}}{\sigma _{\phi }(1,\theta _{T})}\left[
\frac{1}{2n}q_{1,0.05}-KL_{m}(1,\theta _{T})\right] \right)
\label{Approxim}
\end{equation}
where $F_{\mathcal{N}}$ is the cumulative distribution function of
a normal random variable with mean zero and variance one, and
$\phi (x)=-\log x+x-1$; see remarks 3.3 and 3.4 above. The power
function (\ref{Power}) is plotted (with continuous line) for
sample sizes $n=50$, $n=100$, $n=300$ and $n=500$, and for
different values of $\theta_{T}$. Each power entry was obtained
from $1000$ independent runs. The approximation (\ref{Approxim})
is plotted as a function of $\theta_{T}$ by a dashed line. We
observe (see figure 3) that the approximation is accurate for
alternatives which are not ``close to'' the null hypothesis even
for moderate sample sizes.

\begin{figure}[!h]
\centerline{
  \begin{tabular}{ c  c }
 \includegraphics[width=.4\textwidth]{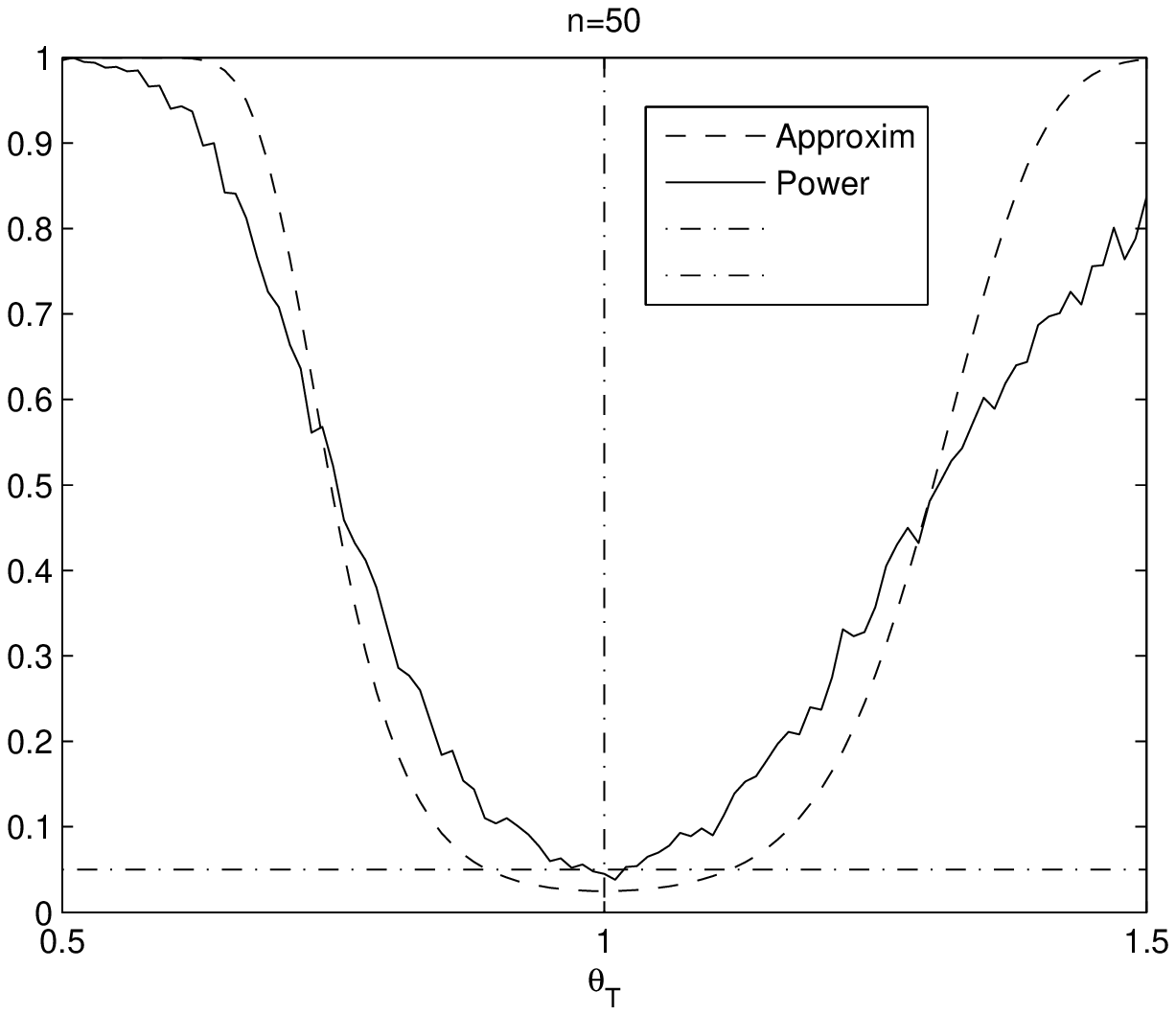}
    &   \includegraphics[width=.4\textwidth]{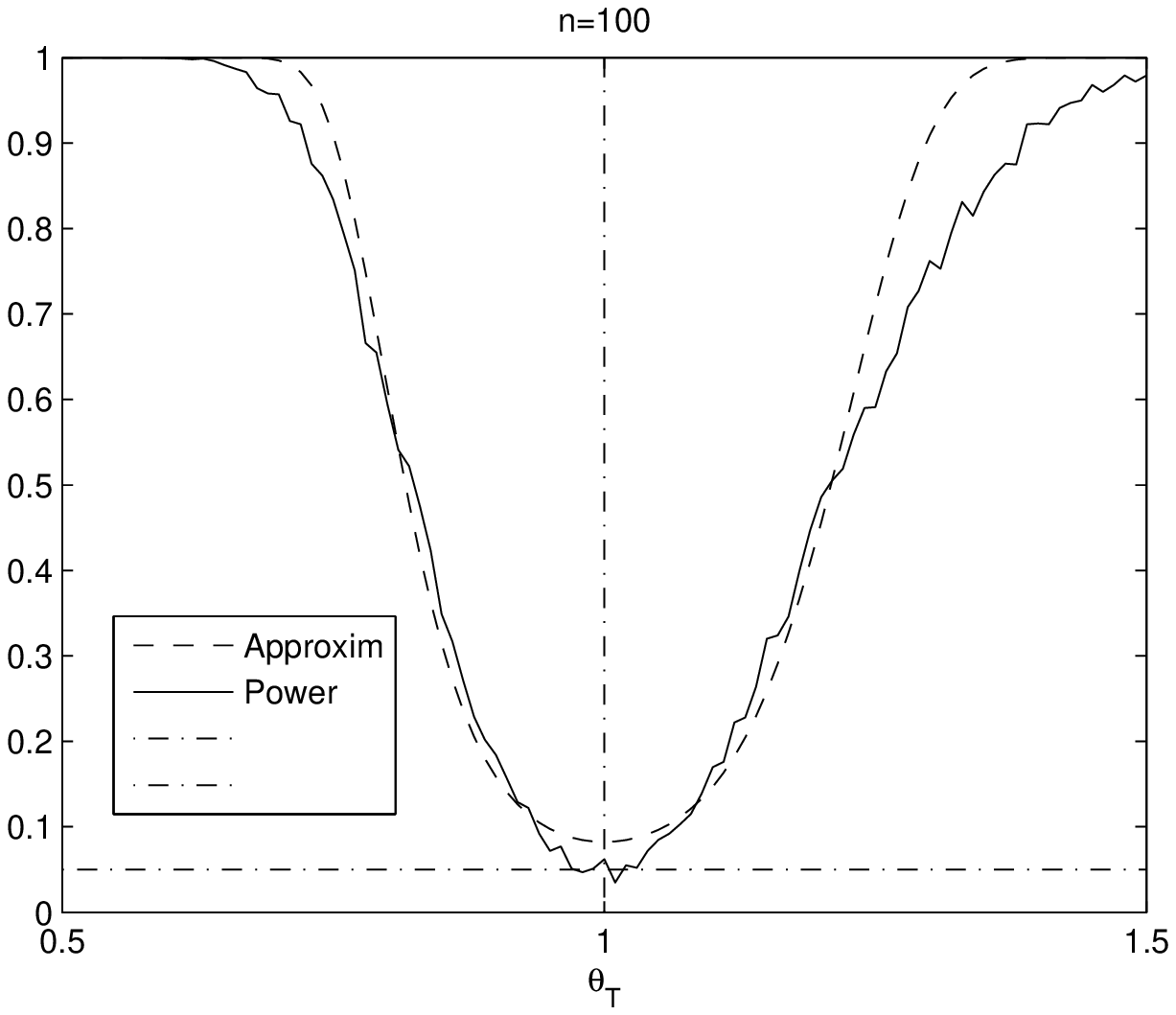}
   \end{tabular} }
\par
\centerline{
  \begin{tabular}{ c  c }
 \includegraphics[width=.4\textwidth]{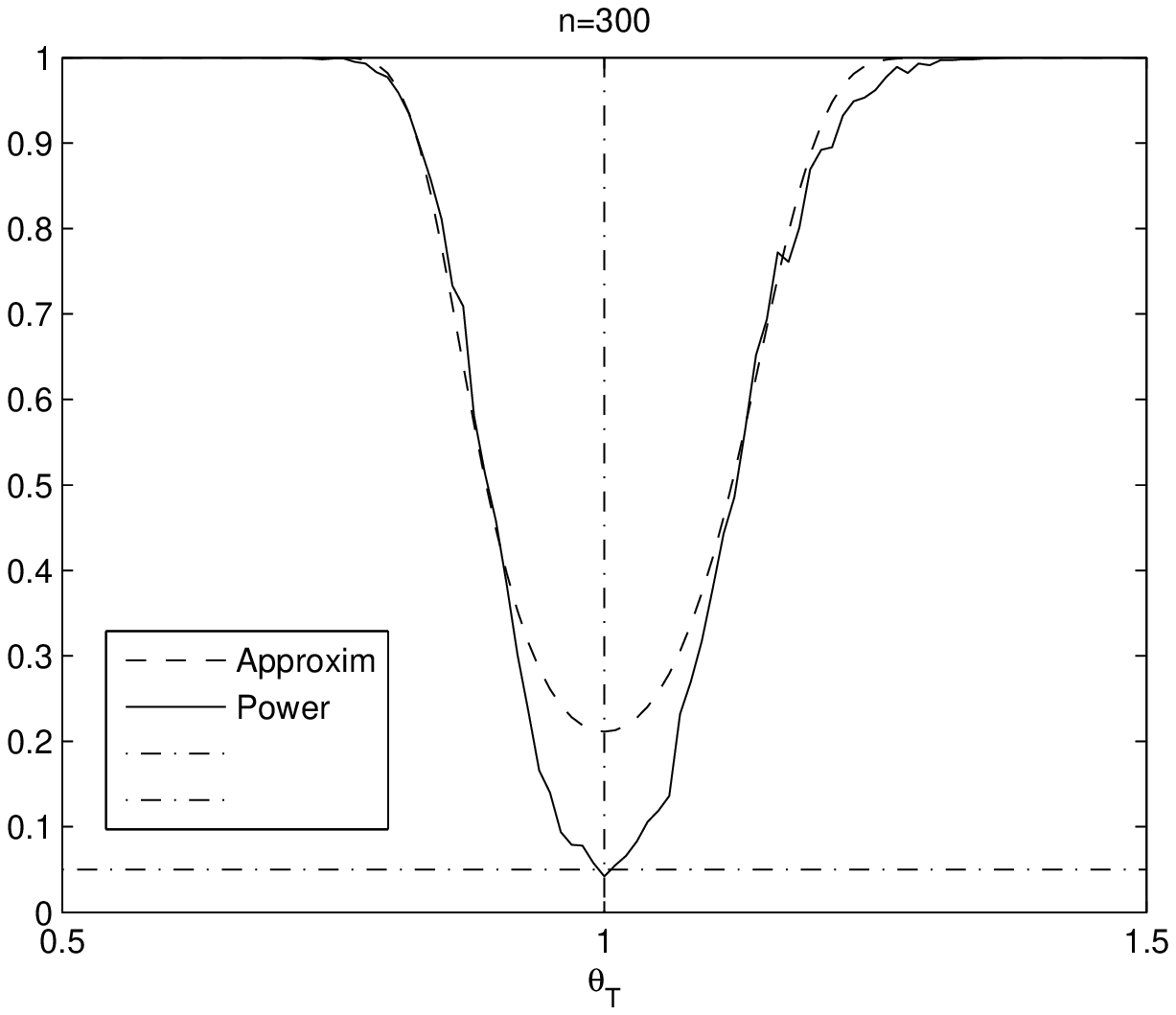}
    &   \includegraphics[width=.4\textwidth]{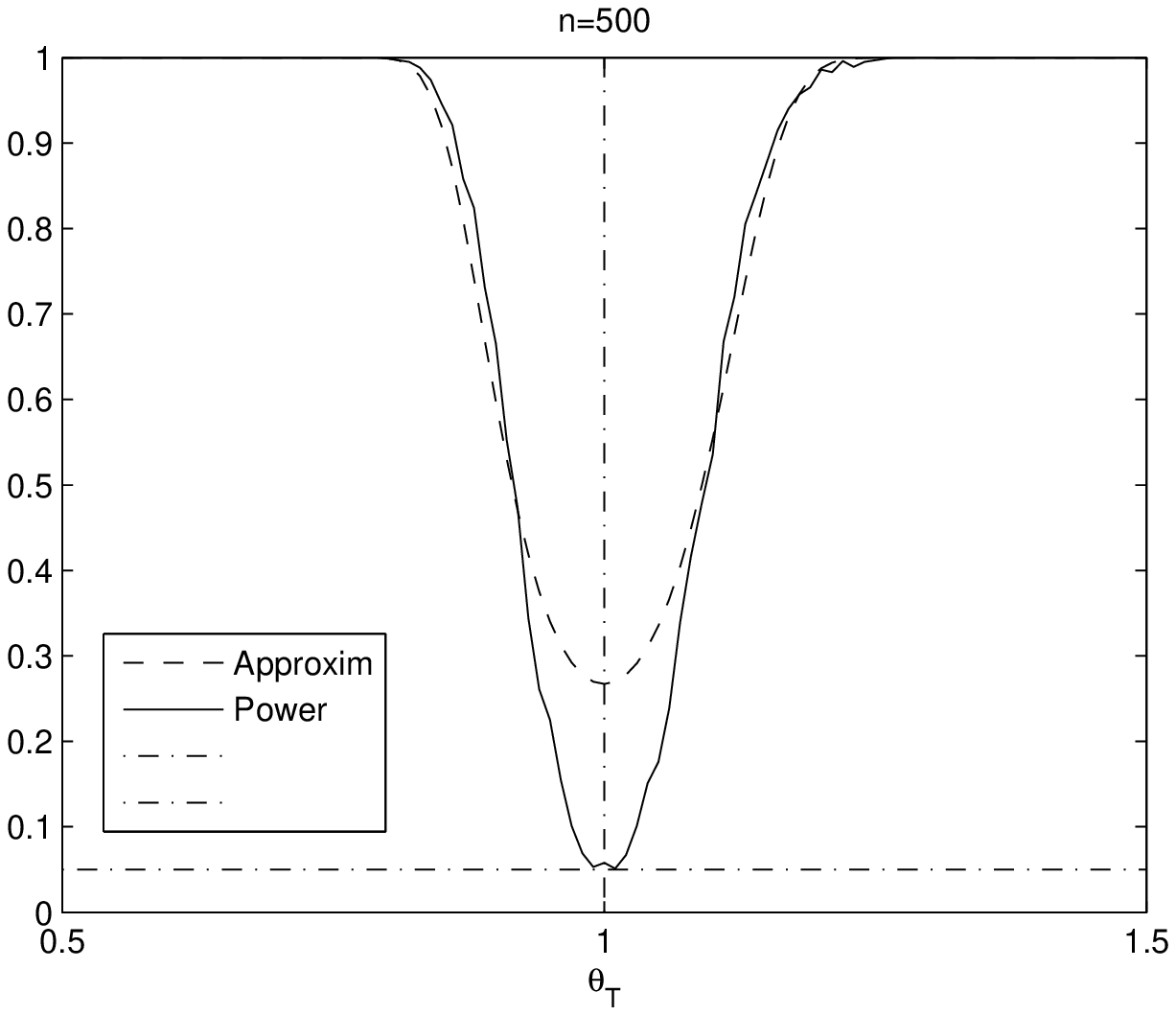}
   \end{tabular} }
\caption{Approximation of the power function}
\label{the limit of the GLRS}
\end{figure}

\section{Concluding remarks and possible developments}

\noindent We have addressed the parametric estimation and test problems. We
have introduced new estimation and test procedure using divergence
minimization and duality technique for discrete or continuous parametric
models, avoiding the smoothing method. The procedure leads to optimal
estimates for the parameter model and for the divergences. It includes both
the discrete (finite or infinite) and the continuous support cases. It
extends the maximum likelihood method for both estimation and test problems.
Moreover, the procedure and the divergences framework permit to obtain the
limit laws of the proposed estimates and the test statistics both under the
null and the alternative (simple or composite) hypotheses, including the
generalized likelihood ratio statistic.  As a by-product, we obtain explicit
power functions in a general case for simple or composite parametric test
problems, and approximations of the minimal sample size which guarantees a
desired power for a given alternative. A new test and new asymptotic
confidence regions are proposed in the case where the parameter may be a
boundary value of the parameter space. Many problems remain to be studied in
the future, such as the choice of the divergence which leads to an
``optimal'' (in some sense) estimate or test in terms of efficiency and
robustness, construction of convergent estimates and test statistics by
divergence when the maximum likelihood is not consistent (for example for
location family for which the expectation does not exists), the Bartlett
correctability and the large deviation properties of the proposed statistics
$\widehat{D_{\phi }}.$

\section{Appendix}

\subsection*{Proof of proposition \ref{prop1}}

(1) We will prove the consistency of the estimate
$\widehat{D_\phi}(\theta,\theta_T)$.  We have
\begin{equation*}
\left|\widehat{D_\phi}(\theta,\theta_T)-D_\phi(\theta,\theta_T)\right|
=\left|P_{n}h(\theta,\widehat{\alpha}_\phi(\theta))-P_{\theta_T}
h(\theta,\theta_T)\right|:=|A|,
\end{equation*}
which implies
\begin{equation*}
P_{n}h(\theta,\theta_T)-P_{\theta_T}h(\theta,\theta_T)\leq A\leq
P_{n}h(\theta,\widehat{\alpha}_\phi(\theta))-
P_{\theta_T}h(\theta,\widehat{\alpha}_\phi(\theta)).
\end{equation*}
Both the RHS and the LHS terms in the above display go to $0$,
under condition (c.2). This implies that $A$ tends to $0$.\\
(2) For the consistency of $\widehat{\alpha}_\phi(\theta)$, we
refer to \cite{vanderVaart1998} theorem 5.7.

\subsection*{Proof of theorem \ref{th asym 1}}
(a) Using (A.1), simple calculus give
\begin{equation}
P_{\theta _{T}}(\partial /\partial \alpha )h(\theta ,\alpha )=0
\label{eqn1}
\end{equation}
and
\begin{equation}
P_{\theta _{T}}(\partial ^{2}/\partial \alpha ^{2})h(\theta
,\theta _{T})=-\int \phi ^{\prime \prime }(p_{\theta }/p_{\theta
_{T}})(p_{\theta }^{2}/p_{\theta _{T}}^{3})p_{\theta _{T}}^{\prime
}p_{\theta _{T}}^{\prime }{}^{T}~d\lambda =:-S.  \label{eqn2}
\end{equation}
Observe that the matrix $S$ is symmetric and positive since the
second derivative $\phi ^{\prime \prime }$ is nonnegative by the
convexity of $\phi $. Let $U_{n}(\theta _{T}):=P_{n}(\partial
/\partial \alpha )h(\theta ,\theta _{T})$, and use (\ref{eqn1})
and (A.2) in connection with the Central Limit Theorem (CLT) to
see that
\begin{equation}
\sqrt{n}U_{n}(\theta _{T})\rightarrow \mathcal{N}(0,M).
\label{eqn3}
\end{equation}%
Also, let $V_{n}(\theta _{T}):=P_{n}(\partial ^{2}/\partial \alpha
^{2})h(\theta ,\theta _{T})$, and use (\ref{eqn2}) and (A.2) in
connection with the Law of Large Numbers (LLN) to conclude that
\begin{equation}
V_{n}(\theta _{T})\rightarrow -S~~(a.s).  \label{eqn4}
\end{equation}%

\noindent  Using the fact that
$P_n(\partial/\partial\alpha)h(\theta,\widehat{\alpha})=0$ and a
Taylor expansion of
$P_n(\partial/\partial\alpha)h(\theta,\widehat{\alpha})$ in
$\widehat{\alpha}$ around $\theta_T$, we obtain
\begin{equation*}
0=P_n(\partial/\partial\alpha)h(\theta,\widehat{\alpha})=
P_n(\partial/\partial\alpha)h(\theta,\theta_T)+(\widehat{\alpha}-\theta_T)^T
P_n(\partial^2/\partial\alpha^2)h(\theta,\theta_T)+o_p(n^{-1/2}).
\end{equation*}
Hence,
\begin{equation}  \label{eqn6}
\sqrt{n}\left(\widehat{\alpha}-\theta_T\right)=
-V_n(\theta_T)^{-1}\sqrt{n} U_n(\theta_T)+o_p(1).
\end{equation}
Using (\ref{eqn3}) and (\ref{eqn4}) and Slutsky theorem, we
conclude then
\begin{equation}
\sqrt{n}\left(\widehat{\alpha}-\theta_T\right)\to \mathcal{N}
\left(0,V_\phi(\theta,\theta_T)\right)
\end{equation}
where $V_\phi(\theta,\theta_T)$ is given in part (a) of theorem
\ref{th asym 1}. When $\theta_T=\theta$, direct calculus shows
that $V_\phi(\theta,\theta_T)=I_{\theta_T}^{-1}$.\newline

\noindent (b) Assume that $\theta_T=\theta$. From (\ref{eqn6}),
using the convergence (\ref{eqn4}), we get
\begin{equation}  \label{eqn7}
\sqrt{n}\left(\widehat{\alpha}-\theta_T\right)=S^{-1}\sqrt{n}U_n(\theta_T)
+o_p(1).
\end{equation}
On the other hand, a Taylor expansion of
$\left[2n/\phi^{\prime\prime}(1)
\right]\widehat{D_\phi}(\theta,\theta_T)=\left[2n/\phi^{\prime\prime}(1)
\right] P_n(\partial/\partial\alpha)h(\theta,\widehat{\alpha})$ in
$\widehat{\alpha}$ around $\theta_T$, using the fact that
$P_nh(\theta,\theta_T)=0$ when $\theta_T=\theta$, gives
\begin{equation*}
\frac{2n}{\phi^{\prime\prime}(1)}\widehat{D_\phi}(\theta,\theta_T)=\frac{2n}{\phi^{\prime\prime}(1)}
U_n^T(\widehat{\alpha}-\theta_T)+\frac{2n}{\phi^{\prime\prime}(1)}(\widehat{\alpha}-\theta_T)^T
V_n(\widehat{\alpha} -\theta_T)+o_p(1).
\end{equation*}
Use (\ref{eqn4}), (\ref{eqn7}) and the fact that
$S=-\phi^{\prime\prime}(1)I_{\theta_T}$ when $\theta_T=\theta$ to
conclude that
\begin{equation*}
\frac{2n}{\phi^{\prime\prime}(1)}\widehat{D_\phi}(\theta,\theta_T)={\phi^{\prime\prime}(1)}^{-2}
\sqrt{n}U_n^TI_{\theta_T}^{-1}\sqrt{n} U_n+o_p(1).
\end{equation*}
Finally, use the convergence (\ref{eqn3}) and the fact that
$M=\phi^{\prime\prime}(1)^2I_{\theta_T}$ when $\theta=\theta_T$,
to conclude that
$\left[2n/\phi^{\prime\prime}(1)\right]\widehat{D_\phi}(\theta,\theta_T)
$ converges in distribution to a $\chi^2$ variable with $d$
degrees of freedom when $\theta=\theta_T$.\newline

\noindent (c) Assume that $\theta_T\neq \theta$. A Taylor
expansion of
$\widehat{D_\phi}(\theta,\theta_T)=P_nh(\theta,\widehat{\alpha})$,
in $\widehat{\alpha}$ around $\theta_T$, using the fact that
$P_{\theta_T}(\partial/\partial\alpha)h(\theta,\theta_T)=0$, gives
$\widehat{D_\phi}(\theta,\theta_T)=P_nh(\theta,\theta_T)+o_p(n^{-1/2})$.
Hence,
\begin{equation*}
\sqrt{n}\left(\widehat{D_\phi}(\theta,\theta_T)-D_\phi(\theta,\theta_T)\right)
=
\sqrt{n}\left[P_nh(\theta,\theta_T)-P_{\theta_T}h(\theta,\theta_T)
\right]+o_p(1),
\end{equation*}
which under assumption (A.3), by the CLT, converges in
distribution to a centred normal variable with variance
$\sigma_\phi^2(\theta,\theta_T)=P_{
\theta_T}h(\theta,\theta_T)^2-\left(P_{\theta_T}
h(\theta,\theta_T)\right)^2$.


\subsection*{Proof of theorem \ref{th asym 1-2}}
(a) For any $\alpha =\theta _{T}+un^{-1/3}$ with $|u|\leq 1$,
consider a Taylor expansion of $P_{n}h(\theta ,\alpha )$ in
$\alpha $ around $\theta _{T}$, and use (A.1) to see that
\begin{equation*}
nP_{n}h(\theta ,\alpha )-nP_{n}h(\theta ,\theta
_{T})=n^{2/3}u^{T}U_{n}+2^{-1}n^{1/3}u^{T}V_{n}u+O(1)~(a.s.)
\end{equation*}%
uniformly on $u$ with $|u|\leq 1$. Now, use (\ref{eqn4}) and the fact that $%
U_{n}=O\left( n^{-1/2}(\log \log n)^{1/2}\right) $ (a.s) to conclude that
\begin{equation*}
nP_{n}h(\theta ,\alpha )-nP_{n}h(\theta ,\theta _{T})=O\left( n^{1/6}(\log
\log n)^{1/2}\right) -2^{-1}u^{T}Sun^{1/3}+O(1)~(a.s.)
\end{equation*}%
uniformly on $u$ with $|u|\leq 1$. Hence, uniformly on the surface of the
ball $B$ (i.e., uniformly on $u$ with $|u|=1$), we have
\begin{equation}
nP_{n}h(\theta ,\alpha )-nP_{n}h(\theta ,\theta _{T})\leq O\left(
n^{1/6}(\log \log n)^{1/2}\right) -2^{-1}cn^{1/3}+O(1)~~(a.s.)  \label{eqn5}
\end{equation}%
where $c$ is the smallest eigenvalue of the matrix $S$. Note that
$c$ is positive since $S$ is positive definite (it is symmetric,
positive and non singular by assumption A.2). In view of
$(\ref{eqn5})$, by the continuity of $\alpha \mapsto P_{n}h(\theta
,\alpha )-nP_{n}h(\theta ,\theta _{T})$ and since it takes value
zero on $\alpha =\theta _{T}$ and is asymptotically negative on
the surface of $B$, it holds that as $n\rightarrow \infty $, with
probability one, $\alpha \mapsto P_{n}h(\theta ,\alpha )$ attains
its maximum value at some point $\widetilde{\alpha}_{\phi }(\theta
)$ in the interior of the ball $B$, and therefore the estimate
$\widetilde{\alpha }_{\phi }(\theta )$ satisfies $P_{n}(\partial
/\partial \alpha )h(\theta ,\widetilde{\alpha })=0$ and
$\widetilde{\alpha }-\theta _{T}=O(n^{-1/3})$.

\noindent The proofs of parts (b), (c) and (d) are similar to
those of parts (a), (b)  and (d) in theorem \ref{th asym 1}.
Hence, they are omitted.

\subsection*{Proof of proposition \ref{prop2}}
We prove (1). For all $\theta \in \Theta $, under condition
(c.4-5-6), we prove that $\sup_{\theta \in \Theta }\Vert
\widehat{\alpha }_{\phi }(\theta )-\theta _{T}\Vert $ tends to
$0$. By the very definition of $\widehat{\alpha }_{\phi }(\theta
)$ and the condition (c.5), we have
\begin{eqnarray}
P_{n}h(\theta ,\widehat{\alpha }_{\phi }(\theta )) &\geq &P_{n}h(\theta
,\theta _{T})  \notag \\
&\geq &P_{\theta _{T}}h(\theta ,\theta _{T})-o_{p}(1),  \notag
\end{eqnarray}
where $o_{p}(1)$ does not depend upon $\theta $ (due to condition
(c.5)). Hence, we have for all $\theta \in \Theta $
\begin{equation}\label{majoration}
P_{\theta _{T}}h(\theta ,\theta _{T})-P_{\theta _{T}}h(\theta
,\widehat{\alpha }_{\phi }(\theta ))\leq P_{n}h(\theta
,\widehat{\alpha }_{\phi }(\theta ))-P_{\theta _{T}}h(\theta
,\widehat{\alpha }_{\phi }(\theta ))+o_{p}(1).
\end{equation}
The  RHS term is less than $\sup_{\left\{ \theta ,\alpha \in
\Theta \right\} }\left| P_{n}h(\theta ,\alpha )-P_{\theta
_{T}}h(\theta ,\alpha )\right| $ $+o_{p}(1)$ which, by (c.5),
tends to $0$. Let $\epsilon >0$ be such that $\sup_{\theta \in
\Theta }\Vert \widehat{\alpha }_{\phi }(\theta )-\theta _{T}\Vert
>\epsilon $. There exists some $a_{n}\in \Theta $ such that $\Vert
\widehat{\alpha }_{\phi }(a_{n})-\theta _{T}\Vert >\epsilon $.
Together with $(c.5.a)$, there exists some $\eta >0$ such that
$P_{\theta _{T}}h(a_{n},\theta _{T})-P_{\theta
_{T}}h(a_{n},\widehat{\alpha }_{\phi }(a_{n}))>\eta $. We then
conclude that
\begin{equation*}
P\left\{ \sup_{\theta \in \Theta }\Vert \widehat{\alpha }_{\phi }(\theta
)-\theta _{T}\Vert >\epsilon \right\} \leq P\left\{ P_{\theta
_{T}}h(a_{n},\theta _{T})-P_{\theta _{T}}h(a_{n},\widehat{\alpha }_{\phi
}(\theta ))>\eta \right\} ,
\end{equation*}
and the RHS term tends to $0$ by (\ref{majoration}). This
concludes the proof of part (1).\newline We prove (2). By the very
definition of $\widehat{\theta }_{\phi }$, conditions  (c.5) and
(c.6) and part (1), we have
\begin{eqnarray}
P_{n}h(\widehat{\theta }_{\phi },\widehat{\alpha }_{\phi
}(\widehat{\theta}_{\phi })) &\leq &P_{n}h(\theta
_{T},\widehat{\alpha }_{\phi }(\theta _{T}))
\notag \\
&\leq &P_{\theta_{T}}h(\theta _{T},\widehat{\alpha }_{\phi
}(\widehat{\theta}_{\phi}))-o_{p}(1),  \notag
\end{eqnarray}%
from which
\begin{eqnarray}\label{maj 2}
P_{\theta _{T}}h(\widehat{\theta }_{\phi },\widehat{\alpha }_{\phi
}(\widehat{\theta }_{\phi }))-P_{\theta _{T}}h(\theta
_{T},\widehat{\alpha }_{\phi }(\widehat{\theta }_{\phi })) &\leq
&P_{\theta _{T}}h(\widehat{\theta }_{\phi },\widehat{\alpha
}_{\phi }(\widehat{\theta }_{\phi
}))-P_{n}h(\widehat{\theta}_{\phi},\widehat{\alpha }_{\phi
}(\widehat{\theta}_{\phi}))+o_{p}(1) \nonumber\\
&\leq &\sup_{\left\{ \theta ,\alpha \in \Theta \right\}}\left|
P_{n}h(\theta ,\alpha )-P_{\theta _{T}}h(\theta,\alpha)\right|
+o_{p}(1).
\end{eqnarray}
Further, by part (1) and condition (c.5.b), for any positive
$\epsilon$, there exists $\eta >0$ such that
\begin{equation*}
P\left\{ \Vert \widehat{\theta }_{\phi }-\theta _{T}\Vert
>\epsilon \right\} \leq P\left\{ P_{\theta _{T}}h(\widehat{\theta
}_{\phi },\widehat{\alpha} _{\phi
}(\widehat{\theta}_{\phi}))-P_{\theta _{T}}h(\theta
_{T},\widehat{\alpha}_{\phi}(\widehat{\theta}_{\phi}))>\eta
\right\},
\end{equation*}
and the RHS term, under condition (c.5), tends to $0$ by (\ref{maj
2}). This concludes the proof.


\subsection*{Proof of theorem \ref{th asym 2}}
Under condition (A.5), simple calculus give
\begin{equation}
P_{\theta _{T}}\frac{\partial }{\partial \alpha }h(\theta
_{T},\theta _{T})=P_{\theta _{T}}\frac{\partial }{\partial \theta
}h(\theta _{T},\theta _{T})=P_{\theta _{T}}\frac{\partial
^{2}}{\partial \alpha \partial \theta} h(\theta _{T},\theta
_{T})=P_{\theta _{T}}\frac{\partial ^{2}}{\partial \theta \partial
\alpha }h(\theta _{T},\theta _{T})=0,  \label{eqn1 1}
\end{equation}
\begin{equation}
-P_{\theta _{T}}\frac{\partial ^{2}}{\partial \alpha ^{2}}h(\theta
_{T},\theta _{T})=P_{\theta _{T}}\frac{\partial ^{2}}{\partial
\theta ^{2}} h(\theta _{T},\theta _{T})=\phi ^{\prime \prime
}(1)I_{\theta _{T}}, \label{eqn1 2}
\end{equation}
and
\begin{eqnarray}
P_{\theta _{T}}\left[ \frac{\partial }{\partial \theta }h(\theta
_{T},\theta _{T})\right] \left[ \frac{\partial }{\partial \theta
}h(\theta _{T},\theta _{T})\right] ^{T} &=&P_{\theta _{T}}\left[
\frac{\partial }{\partial \alpha } h(\theta _{T},\theta
_{T})\right] \left[ \frac{\partial }{\partial \alpha }
h(\theta _{T},\theta _{T})\right] ^{T}  \notag \\
&=&-P_{\theta _{T}}\left[ \frac{\partial }{\partial \alpha
}h(\theta _{T},\theta _{T})\right] \left[ \frac{\partial
}{\partial \theta }h(\theta
_{T},\theta _{T})\right] ^{T}  \notag \\
&=&{\phi ^{\prime \prime }(1)}^{2}I_{\theta _{T}}.
\end{eqnarray}
Denote $U_{n}(\theta ,\theta _{T}):=P_{n}(\partial /\partial
\alpha )h(\theta ,\theta _{T})$, $V_{n}(\theta ,\theta
_{T}):=P_{n}(\partial ^{2}/\partial \alpha ^{2})h(\theta ,\theta
_{T})$, $S(\theta ,\theta _{T}):=-P_{\theta _{T}}(\partial
^{2}/\partial \alpha ^{2})h(\theta ,\theta _{T})$ and
$a_{n}^{T}:=\left( (\widehat{\theta }_{\phi }-\theta
_{T})^{T},(\widehat{\alpha }_{\phi }(\widehat{\theta }_{\phi
})-\theta _{T})^{T}\right) ^{T}$. Under conditions (A.4-5), by a
Taylor expansion, we obtain
\begin{equation*}
\sqrt{n}a_{n}=\sqrt{n}\left[
\begin{array}{cc}
\frac{1}{\phi ^{\prime \prime }(1)}I_{\theta _{T}}^{-1} & 0 \\
0 & \frac{-1}{\phi ^{\prime \prime }(1)}I_{\theta _{T}}^{-1}
\end{array}
\right] \left[
\begin{array}{c}
-P_{n}\frac{\partial }{\partial \theta }h(\theta _{T},\theta _{T}) \\
-P_{n}\frac{\partial }{\partial \alpha }h(\theta _{T},\theta _{T})
\end{array}
\right] +o_{p}(1).
\end{equation*}
We therefore deduce, by the CLT, that, under condition (A.6),
$\sqrt{n}a_{n}$ converges in distribution to a centred normal
variable with covariance matrix
\begin{equation*}
\mathbb{V}=\left[
\begin{array}{cc}
I_{\theta _{T}}^{-1} & I_{\theta _{T}}^{-1} \\
I_{\theta _{T}}^{-1} & I_{\theta _{T}}^{-1}
\end{array}
\right],
\end{equation*}
which completes the proof of theorem \ref{th asym 2}.


\subsection*{Proof of theorem \ref{th asym 2-2}}
(a) Using condition (A.5) and (\ref{eqn1 1}), we can write
\begin{equation}
U_{n}(\theta ,\theta _{T}):=U_{n}(\theta _{T},\theta
_{T})+o(n^{-1/3})~(a.s.) \label{eqn2 1}
\end{equation}
and
\begin{equation}
V_{n}(\theta ,\theta _{T}):=V_{n}(\theta _{T},\theta
_{T})+O(n^{-1/3})~(a.s.),  \label{eqn2 2}
\end{equation}
uniformly on $\theta \in B(\theta _{T},n^{-1/3})$. On the other
hand, for any $\alpha =\theta _{T}+un^{-1/3}$ with $|u|\leq 1$, by
a Taylor expansion using condition (A.5), we obtain
\begin{equation*}
nP_{n}h(\theta ,\alpha )-nP_{n}h(\theta ,\theta
_{T})=n^{2/3}u^{T}U_{n}(\theta ,\theta
_{T})+2^{-1}n^{1/3}u^{T}V_{n}(\theta ,\theta _{T})u+O(1)~(a.s.)
\end{equation*}
uniformly on $\theta \in B(\theta _{T},n^{-1/3})$ and $u$ with $|u|\leq 1$.
Combining this with (\ref{eqn2 1}) and (\ref{eqn2 2}) to see that
\begin{equation*}
nP_{n}h(\theta ,\alpha )-nP_{n}h(\theta ,\theta
_{T})=n^{2/3}u^{T}U_{n}(\theta _{T},\theta
_{T})+2^{-1}n^{1/3}u^{T}V_{n}(\theta _{T},\theta
_{T})u+o(n^{1/3})~(a.s.)
\end{equation*}
uniformly on $\theta \in B(\theta _{T},n^{-1/3})$ and $u$ with $|u|\leq 1$.
Now, from this, using the fact that $U_{n}(\theta _{T},\theta _{T})=O\left(
n^{-1/2}(\log \log n)^{1/2}\right) $ (a.s.) and $V_{n}(\theta _{T},\theta
_{T})=-S(\theta _{T},\theta _{T})+o(1)$ (a.s.), we obtain
\begin{equation}
nP_{n}h(\theta ,\alpha )-nP_{n}h(\theta ,\theta _{T})=O\left(
n^{1/6}(\log \log n)^{1/2}\right) -2^{-1}n^{1/3}u^{T}S(\theta
_{T},\theta _{T})u+o(n^{1/3})~(a.s.)  \label{eqn2 3}
\end{equation}
uniformly on $\theta \in B(\theta _{T},n^{-1/3})$ and $u$ with $|u|\leq 1$.
Hence, uniformly on $\alpha $ in the surface of the ball $B(\theta
_{T},n^{-1/3})$ (i.e., uniformly on $u$ with $|u|=1$), we have
\begin{equation}
nP_{n}h(\theta ,\alpha )-nP_{n}h(\theta ,\theta _{T})\leq O\left(
n^{1/6}(\log \log n)^{1/2}\right) -2^{-1}\phi ^{\prime \prime
}(1)cn^{1/3}+o(n^{1/3})~(a.s.)  \label{eqn2 3.1}
\end{equation}
(uniformly on $\theta \in B(\theta _{T},n^{-1/3}))$ where $c>0$ is
the smallest eigenvalue of the matrix $I_{\theta _{T}}=\phi
^{\prime \prime }(1)^{-1}S(\theta _{T},\theta _{T})$. Hence, by
the continuity of the function $\alpha\mapsto nP_{n}h(\theta
,\alpha )-nP_{n}h(\theta ,\theta _{T})$ and since it takes value
zero when  $\alpha=\theta_T$ and is asymptotically negative with
respect to $\alpha$ on  the surface of $B$, it holds that, as $n$
tends to $\infty$, with probability one,  the function
$\alpha\mapsto P_nh(\theta,\alpha)$ attains it maximum value at
some point $\widetilde{\alpha}_\phi(\theta)$ in the interior of
$B(\theta_T,n^{-1/3})$, and this holds for all $\theta\in
B(\theta_T,n^{-1/3})$. Further, since (\ref{eqn2 3}) holds
uniformly on $\theta\in B(\theta_T,n^{-1/3})$, we conclude that
\begin{equation}\label{eqn2 4}
\widetilde{\alpha}_\phi(\theta)-\theta_T=O(n^{-1/3}) ~(a.s.)
\text{ uniformly on } \theta\in B(\theta_T,n^{-1/3}).
\end{equation}
\noindent We now prove that, as $n\rightarrow \infty $, with
probability one, the function $\theta \mapsto P_{n}(\theta
,\widehat{\alpha }_{\phi }(\theta ))$ attains its minimum value at
some point $\widetilde{\theta }_{\phi }$ in the interior of the
ball $B(\theta _{T},n^{-1/3})$. Here, $\widetilde{\alpha} _{\phi
}(\theta )$ is any value in the interior of $B(\theta
_{T},n^{-1/3})$ which maximizes $\alpha \mapsto P_{n}h(\theta
,\alpha )$. It exists by the above arguments. For any $\theta
=\theta _{T}+vn^{-1/3}$ with $|v|\leq 1$, by a Taylor expansion of
$nP_{n}h(\theta ,\widetilde{\alpha }_{\phi }(\theta )) $ in
$\theta $ and $\widetilde{\alpha }_{\phi }(\theta )$ around
$\theta _{T}$ , and a Taylor expansion of $nP_{n}h(\theta
_{T},\widetilde{\alpha }_{\phi }(\theta _{T}))$ in
$\widetilde{\alpha }_{\phi }(\theta _{T})$ around $\theta _{T}$,
using (\ref{eqn2 4}) and (\ref{eqn1 1}), we obtain
\begin{eqnarray*}
nP_{n}h(\theta ,\widetilde{\alpha }_{\phi }(\theta
))-nP_{n}h(\theta _{T}, \widetilde{\alpha }_{\phi }(\theta _{T}))
&=&n^{2/3}v^{T}P_{n}(\partial
/\partial \theta )h(\theta _{T},\theta _{T})+ \\
&&2^{-1}n^{1/3}v^{T}\left[ P_{n}(\partial ^{2}/\partial \theta ^{2})h(\theta
_{T},\theta _{T})\right] v+o(n^{1/3})~(a.s.)
\end{eqnarray*}
uniformly on $v$ with $|v|\leq 1$. Hence, from this, using the
fact that
\newline
$P_{n}(\partial /\partial \theta )h(\theta _{T},\theta _{T})=O\left(
n^{-1/2}(\log \log n)^{1/2}\right) $ (a.s.) and $P_{n}(\partial
^{2}/\partial \theta ^{2})h(\theta _{T},\theta _{T})=\phi ^{\prime \prime
}(1)I_{\theta _{T}}+o(1)$ (a.s.), we conclude that
\begin{equation*}
nP_{n}h(\theta ,\widetilde{\alpha }_{\phi }(\theta
))-nP_{n}h(\theta _{T}, \widetilde{\alpha }_{\phi }(\theta
_{T}))=O\left( n^{1/6}(\log \log n)^{1/2}\right) +2^{-1}\phi
^{\prime \prime }(1)v^{T}I_{\theta _{T}}vn^{1/3}+o(n^{1/3})~(a.s.)
\end{equation*}
uniformly on $v$ with $|v|\leq 1$. Hence, uniformly on $\theta $
in the surface of the ball $B(\theta _{T},n^{-1/3})$ (i.e.,
uniformly on $v$ with $|v|=1$), we obtain
\begin{equation*}
nP_{n}h(\theta ,\widetilde{\alpha }_{\phi }(\theta
))-nP_{n}h(\theta _{T}, \widetilde{\alpha }_{\phi }(\theta
_{T}))\geq O\left( n^{1/6}(\log \log n)^{1/2}\right) +2^{-1}\phi
^{\prime \prime }(1)cn^{1/3}+o(n^{1/3})~(a.s.)
\end{equation*}
where $c>0$ is the smallest eigenvalue of $I_{\theta _{T}}$. This implies
that
\begin{equation*}
n^{2/3}P_{n}h(\theta ,\widetilde{\alpha }_{\phi }(\theta
))-n^{2/3}P_{n}h(\theta _{T},\widetilde{\alpha }_{\phi }(\theta
_{T}))\geq O\left( n^{-1/6}(\log \log n)^{1/2}\right) +2^{-1}\phi
^{\prime \prime }(1)c+o(1)~(a.s.)
\end{equation*}
uniformly on $\theta $ in the surface of the ball $B(\theta
_{T},n^{-1/3})$. The left hand side of the above display equals
zero when $\theta =\theta _{T} $ and is positive when $\theta $ is
in the surface of the ball $B(\theta _{T},n^{-1/3})$ (for $n$
sufficiently large). This implies that, as $n\rightarrow \infty$,
with probability one, the function $\theta \mapsto P_{n}h(\theta
,\widetilde{\alpha }_{\phi }(\theta ))$ attains its minimum value
at some point $\widetilde{\theta }_{\phi }$ in the interior of the
ball $B$. This concludes the proof of part (a).

\noindent (b) See the proof of theorem \ref{th asym 2}.

\subsection*{Proof of theorem \ref{th asym 3}}
We have
\begin{eqnarray*}
\widehat{D_\phi}(\Theta_0,\theta_T) & := & \inf_{\beta \in
B_{0}}\sup_{\alpha \in \Theta }P_{n}h\left(s(\beta),\alpha\right). \\
& = & P_{n}h\left(s(\widehat{\beta}),\widehat{\alpha}\right),
\end{eqnarray*}
in which as in the proof of theorem \ref{th asym 2}, $s(\widehat{\beta})$
and $\widehat{\alpha}$ are solutions of the system of equations
\begin{equation*}
\left\{
\begin{array}{ccc}
P_{n}\frac{\partial }{\partial \beta
}h\left(s(\widehat{\beta}),\widehat{\alpha}\right) & = & 0 \\
P_{n}\frac{\partial }{\partial \alpha }h\left(
s(\widehat{\beta}),\widehat{\alpha}\right) & = & 0.
\end{array}
\right.
\end{equation*}
In the first equation the partial derivative is intended w.r.t. the first
variable $\beta$ in $s(\beta)$ and in the second one w.r.t. the second
variable $\alpha$. A Taylor expansion of $P_{n}\frac{\partial }{\partial
\beta }h\left(s(\widehat{\beta}),\widehat{\alpha}\right)$ and $P_{n}\frac{%
\partial }{\partial \alpha }h\left( s(\widehat{\beta}),\widehat{\alpha}%
\right)$ in a neighborhood of $(\beta_T,\theta_T)$ gives
\begin{equation}  \label{eqn theta beta}
\left[
\begin{array}{c}
-P_{n}\frac{\partial }{\partial \beta }h(s(\beta_T),\theta_T) \\
-P_{n}\frac{\partial }{\partial \alpha }h(s(\beta_T),\theta_T)
\end{array}
\right] =\left[
\begin{array}{cc}
P_{\theta_T}\frac{\partial
^{2}}{\partial\beta^2}h(s(\beta_T),\theta_T) &
P_{\theta_T}\frac{\partial ^{2}}{\partial \beta \partial \alpha }
h(s(\beta_T),\theta_T) \\
P_{\theta_T}\frac{\partial ^{2}}{\partial \alpha \partial \beta }
h(s(\beta_T),\theta_T) & P_{\theta_T}\frac{\partial ^{2}}{\partial
\beta^{2}} h(s(\beta_T),\theta_T)
\end{array}
\right]b_n+o_p(1),
\end{equation}
where $b_n:=\left((\widehat{\beta}-\beta_T)^T,(\widehat{\alpha}
-\theta_T)^T\right)^T$. This implies that $b_n=O_p(n^{-1/2})$. So,
by a Taylor expansion of $\widehat{D_\phi}(\Theta_0,\theta_T)$
around $(\beta_T,\theta_T)$, we obtain
\begin{equation}
\frac{2n}{\phi ^{\prime \prime }(1)}T_{n}^{\phi
}=U_{n}^{T}A^{-1}U_{n}-V_{n}^{T}B^{-1}V_{n}+o_p(1),  \label{stat
composite}
\end{equation}
where
\begin{equation*}
U_{n} := \frac{\sqrt{n}}{\phi ^{\prime \prime }(1)}P_{n}\frac{\partial }{
\partial \alpha }h\left(s(\beta_T),\theta_T\right),\quad V_{n} :=
\frac{\sqrt{n}}{\phi ^{\prime \prime }(1)}P_{n}\frac{\partial }{
\partial \beta} h(s(\beta_T),\theta_T),
\end{equation*}
\begin{equation*}
A := -\frac{1}{\phi ^{\prime \prime
}(1)}P_{\theta_T}\frac{\partial ^{2}}{\partial
\alpha^{2}}h(s(\beta_T),\theta_T),\quad B := \frac{1}{\phi
^{\prime \prime }(1)}P_{\theta _T}\frac{\partial ^{2}}{\partial
\beta ^{2}} h(s(\beta_T),\theta_T).
\end{equation*}
By (\ref{eqn1 2}), it holds $A=I_{\theta_T}$. On the other hand,
\begin{eqnarray}
\frac{\partial }{\partial \beta }h\left(s(\beta_T),\theta_T\right)
& = & {\left[ \frac{\partial }{\partial \beta }s(\beta_T)\right]
}^{T}\frac{\partial }{ \partial s(\beta )}h\left(s(\beta_T),\theta_T\right)  \notag \\
& = & {\left[ S(\beta_T)\right] }^{T}\frac{\partial }{\partial s(\beta)}
h\left(s(\beta_T),\theta_T\right).  \notag
\end{eqnarray}
Moreover, using the fact that $\phi ^{\prime }(1)=0$, we can see
that $\frac{\partial }{\partial
s(\beta)}h\left(s(\beta_T),\theta_T\right)=-\frac{\partial
}{\partial \alpha } h\left(s(\beta_T),\theta_T\right)$, which
implies
\begin{equation*}
P_{\theta_T}\frac{\partial }{\partial \beta
}h\left(s(\beta_T),\theta_T \right)={\left[ S(\beta_T)\right] }
^{T}\left[ -P_{\theta_T}\frac{\partial }{\partial \alpha
}h\left(s(\beta_T),\theta_T\right)\right].
\end{equation*}
In the same way, we obtain
\begin{equation*}
P_{\theta_T}\frac{\partial ^{2}}{\partial \beta
^{2}}h\left(s(\beta_T), \theta_T\right)={\left[ S(\beta_T)\right]
}^{T}\left[ -P_{\theta_T} \frac{\partial ^{2}}{\partial
\alpha^2}h\left(s(\beta_T),\theta_T\right)\right] \left[
S(\beta_T)\right].
\end{equation*}
It follows that $V_{n}={\left[ S(\beta_T)\right] }^{T}U_{n}$ and
$B={\left[ S(\beta_T)\right] }^{T}I_{\theta_T}S(\beta_T)$.
Combining this result with (\ref{stat composite}), we get
\begin{equation*}
\frac{2n}{\phi ^{\prime \prime }(1)}\widehat{D_\phi}(\Theta_0,
\theta_T)=U_{n}^{T}\left[
I_{\theta_T}^{-1}-S(\beta_T)B^{-1}{S(\beta_T)}^{T} \right]
U_{n}+o_p(1),
\end{equation*}
which is precisely the asymptotic expression for the Wilks likelihood ratio
statistic for composite hypotheses. The proof is completed following
therefore the same arguments as for the Wilks likelihood ratio statistic;
see e.g. \cite{Sen-Singer1993} chapter 5.

\subsection*{Proof of theorem \ref{th asym 4}}

 The proofs of part (a) and (b) are similar to the proofs of part (a) and
(b) of theorem \ref{th asym 3}, hence they are omitted.\newline

\noindent (c) Using (\ref{formule de base}) and (\ref{formule de
base simple} ), we can see that $D_\phi(\Theta_0,\theta_T)$ can be
written as
\begin{eqnarray}  \label{eqn2 5}
D_\phi(\Theta_0,\theta_T) & := & \inf_{\beta\in B_0}
D_\phi(s(\beta),\theta_T) = D_\phi(s(\beta^*),\theta_T)  \notag \\
& = & \sup_{\alpha\in\Theta} P_{\theta_T}h(s(\beta^*),\alpha) =
P_{\theta_T}h(s(\beta^*),\theta_T).
\end{eqnarray}
On the other hand, by a Taylor expansion of $\widehat{D_\phi}
(\Theta_0,\theta_T)=P_nh(s(\widehat{\beta}),\widehat{\alpha}_\phi(\widehat{\beta}))$
in $\widehat{\beta}$ and $\widehat{\alpha}_\phi(\widehat{\beta})$
around $\beta^*$ and $\theta_T$, we obtain
\begin{equation*}
\widehat{D_\phi}(\Theta_0,\theta_T)= P_nh(s(\beta^*),\theta_T)+o_p(n^{-1/2}).
\end{equation*}
Combining this with (\ref{eqn2 5}) to conclude that
\begin{equation*}
\sqrt{n}\left[\widehat{D_\phi}(\Theta_0,\theta_T)-D_\phi(\Theta_0,\theta_T)\right]=
\sqrt{n}\left[P_nh(s(\beta^*),\theta_T)-P_{\theta_T}h(s(\beta^*),
\theta_T)\right]+o_p(1)
\end{equation*}
which, by the CLT, converges to a centred normal variable with variance
\begin{equation*}
\sigma_\phi^2(\beta^*,\theta_T)= P_{\theta_T}h(s(\beta^*),\theta_T)^2-
\left(P_{\theta_T}h(s(\beta^*),\theta_T)\right)^2.
\end{equation*}
This ends the proof.

\end{document}